%% file: main.tex
\documentclass[11pt,reqno]{amsart}
\usepackage[T1]{fontenc}
\usepackage{lmodern}
\usepackage{amsmath,amssymb,amsthm,mathtools}
\usepackage{booktabs,array,graphicx,multirow}
\usepackage{flafter}
\usepackage[numbers,square,sort&compress]{natbib}
\usepackage[colorlinks=true,linkcolor=blue,citecolor=blue,urlcolor=blue]{hyperref}
\usepackage{xcolor}
\numberwithin{equation}{section}
\newcommand{\pinnTheoremName}{Theorem}
\newtheorem{theorem}{\pinnTheoremName}[section]

\theoremstyle{remark}
\newtheorem{remark}[theorem]{Remark}
\newcommand{\R}{\mathbb R}
\newcommand{\eps}{\varepsilon}
\newcommand{\dd}{\,\mathrm d}
\newcommand{\norm}[1]{\lVert#1\rVert}

\newcommand{\Lstd}{\mathcal L_{\mathrm{std}}}
\newcommand{\Laug}{\mathcal L_{\mathrm{aug}}}
\newcommand{\E}{\mathcal E}
\title[Integral constraints in physics-informed neural networks]{When and why can integral constraints improve the accuracy of physics-informed neural networks?}
\author{Myeong-Su Lee}
\address{Research Institute of Mathematics, Seoul National University, Seoul, \mbox{Republic of Korea}}
\email{msl3573@snu.ac.kr}
\date{}
\begin{document}
\begin{abstract}
Adding constraints obtained by integrating governing equations to the loss
has been reported to improve the accuracy of physics-informed neural
networks (PINNs) in several applications, but mathematical understanding
of these benefits remains limited. For several partial differential
equations (PDEs), we derive stability estimates for losses with and
without integral constraints. The comparison reveals that in certain
parameter regimes, the standard loss has very weak control over specific
error components, and that the integral constraints restore uniform
control of these components. This explains when and why the constraints
improve accuracy. Finally, numerical experiments support this
interpretation. In parameter regimes where the standard loss is predicted
to have weak control over those components, the two losses
indeed show clear differences in accuracy. Conversely, in regimes where
the standard loss retains sufficient control, both losses yield similar
accuracy.
\end{abstract}
\maketitle
\enlargethispage{3pt}

\section{Introduction}\label{sec:introduction}
Physics-informed neural networks (PINNs) approximate the solution of a
partial differential equation (PDE) with a neural network. The network is
trained by minimizing a loss based on the governing PDE and
given data, such as initial and boundary conditions
\citep{RaissiEtAl2019}. PINNs have been applied to forward
and inverse problems across a wide range of physics and engineering
disciplines, from fluid mechanics to battery modelling
\citep{RaissiYazdaniKarniadakis2020,KissasEtAl2020,ChenEtAl2020,Haghighat2021,ChoTraveling2021,Patel2022,ChangEtAl2024,DickEtAl2025,LeeEtAl2025,ChoConformal2026}.
Alongside these successful applications, PINNs have been actively
studied from various theoretical perspectives
\citep{DeRyckMishra2024,DoumecheConvergence2025,XuLiHuang2025}.

Integrating a PDE and using the given initial and boundary
conditions can give additional relations that its solution must satisfy.
Depending on the PDE, these relations are referred to as conservation
laws, balance relations, or compatibility conditions. Adding penalties for violations of these
relations has been reported to improve PINN accuracy across a range of
applications. However, explanations of the gains are limited and mostly
heuristic
\citep{HennighEtAl2021,LinChen2022,HassanalyEtAl2024,LeeEtAl2025,HuangEtAl2026,ZhangMUSA2026}.
Moreover, in some cases, such penalties have been reported to bring no
improvement or even to reduce accuracy \citep{LinChen2022,BaezEtAl2026}.
Understanding the mechanism behind these observations is therefore
important for the reliable use and further development of such integral
constraints.

Stability theory is a classical topic in the analysis of PDEs.
It studies how sensitive solutions are to perturbations in the PDE or in the initial and boundary conditions.
In the analysis of PINNs, stability estimates have been used to control
the solution error by the residual loss for a range of linear and
nonlinear PDEs
\citep{MishraMolinaro2023,DeRyckJagtapMishra2024,ZeinhoferEtAl2024,DeRyck2024,ErnstUrban2025,GazoulisEtAl2025}.
In this context, for certain PDEs the standard loss has been found
insufficient to control the solution error adequately, leading to
alternative loss formulations tailored to those PDEs
\citep{WangEtAl2022HJB,KoEtAl2026,ChoiKoKwon2026,LuWangXu2022,JinMaWu2023,JinMaWu2024}.
In addition, error estimates have been obtained for losses that
include integral constraints \citep{LuLuWang2021,BaiEtAl2025}. However,
to the best of our
knowledge, the improvement in accuracy has not been explained by a
direct comparison of such estimates with and without the
constraints.

In this work, we use classical stability theory to understand the relation
between the loss and the solution error. From this perspective, we
consider how integral constraints can improve accuracy. For several
elliptic and time-dependent PDEs, we show that in certain parameter
regimes the standard loss has weak control over some error components,
causing the constants in the error estimates to become necessarily large.
This means that reducing the standard loss alone cannot effectively
reduce the solution error in these regimes. On the other hand, adding
integral constraints to the loss restores uniform control of these
components and removes the parameter dependence from the error bounds.
Consequently, when the two losses are reduced to comparable levels, the
loss with the constraints guarantees a much smaller error.

A simple elliptic problem illustrates this mechanism. Let us consider
\begin{equation}\label{eq:intro-neumann}
 -u''+\varepsilon u=f\quad\hbox{in }(0,1),
 \qquad u'(0)=u'(1)=0,\qquad 0<\varepsilon\le1,
\end{equation}
with $f\in L^2(0,1)$. Let $u$ denote its solution. For a trial function $v\in H^2(0,1)$ satisfying
the same boundary conditions, the standard loss is defined by
\begin{equation}\label{eq:intro-standard-loss}
 \mathcal L_{\mathrm{std},\varepsilon}(v)
 :=\|-v''+\varepsilon v-f\|_{L^2(0,1)}^2.
\end{equation}
For the loss \eqref{eq:intro-standard-loss}, the following error bound holds:
\begin{equation}\label{eq:intro-standard}
 \|v-u\|_{L^2(0,1)}^2
 \le\varepsilon^{-2}\mathcal L_{\mathrm{std},\varepsilon}(v).
\end{equation}
This means that for small $\varepsilon$, guaranteeing a target accuracy
requires a very small loss. This weakening of error control arises
because the standard loss becomes insensitive to certain error
components. Specifically, for any constant $c>0$, let us consider
$v:=u+c$. Then, this function satisfies the boundary conditions exactly
and changes only the reaction term in \eqref{eq:intro-neumann}, so the
standard loss equals $\varepsilon^2c^2$. For smaller $\varepsilon$,
training must achieve a much smaller loss to detect such a constant
error.

However, limited network expressivity and
incomplete optimization can prevent training from reaching such a small
loss. In practice, the trained approximation may therefore carry a
constant offset, and the standard loss alone may not guarantee sufficient
accuracy for small $\varepsilon$. To address this, we use the integral
constraint obtained by integrating \eqref{eq:intro-neumann} over the
domain and applying the boundary conditions:
\begin{equation}\label{eq:intro-mean}
 \int_0^1u\,dx=\varepsilon^{-1}\int_0^1f\,dx.
\end{equation}
Measuring how much $v$ violates this condition and adding it
to the standard loss gives the augmented loss
\begin{equation}\label{eq:intro-augmented-loss}
 \mathcal L_{\mathrm{aug},\varepsilon}(v)
 :=\mathcal L_{\mathrm{std},\varepsilon}(v)
 +\left|\int_0^1v\,dx-\varepsilon^{-1}\int_0^1f\,dx\right|^2.
\end{equation}
For the same shift $v=u+c$, the added term is exactly $c^2$,
independently of $\varepsilon$.
In fact, for this augmented loss, the following error bound holds:
\begin{equation}\label{eq:intro-augmented}
 \|v-u\|_{L^2(0,1)}^2
 \le\mathcal L_{\mathrm{aug},\varepsilon}(v).
\end{equation}
The factor $\varepsilon^{-2}$ appearing in
\eqref{eq:intro-standard} has been removed. If training reduces the two
losses to comparable values, \eqref{eq:intro-augmented} provides a much
stronger error bound for small $\varepsilon$. By contrast, when
$\varepsilon$ is close to one, the two bounds are nearly equal, so the
degree of accuracy improvement depends on $\varepsilon$. The numerical
experiments in Section~\ref{sec:numerics} confirm this trend.

Section~\ref{sec:theory} extends this stability analysis to a nonlinear
elliptic equation with Neumann boundary conditions, the Helmholtz and
stationary Burgers equations with Dirichlet boundary conditions, and a
time-dependent reaction--diffusion equation, and
Section~\ref{sec:numerics} presents numerical experiments for these four
PDEs. Our main contributions are as follows.
\begin{itemize}
\item For each PDE, we compare the stability estimates with and without
an integral constraint. In certain parameter regimes, we construct
explicit approximations that attain a small standard loss while retaining
a large error. These constructions show that the standard estimate is
sharp and identify the form of error that the standard loss fails to
detect. The integral constraint controls the corresponding error
component uniformly in the parameter.
\item Numerical experiments for all four equations support this
interpretation. In regimes where the stability estimates predict weak
error control, errors of the forms identified in
Section~\ref{sec:theory} indeed appear in approximations trained with
the standard loss, and training with the augmented loss substantially
reduces them. In regimes where the standard loss already provides
effective control, both losses yield similar accuracy.
\end{itemize}

\section{Preliminaries}\label{sec:prelim}

\subsection{Physics-Informed Neural Networks}\label{sec:pinn-formulation}
We consider a spatial domain $\Omega\subset\mathbb R^d$ and a terminal time $T>0$.
For an evolution problem, we write $Q:=(0,T)\times\Omega$ and
$\Sigma:=(0,T)\times\partial\Omega$. Let $\mathcal F$ and $\mathcal B$ denote the differential and boundary operators, respectively. For given source $f$, boundary data $g$ and initial data $u_0$, the problem takes the form
\begin{equation}\label{eq:pinn-problem}
 \mathcal F(u)=f\quad\text{in }Q,\qquad
 \mathcal B(u)=g\quad\text{on }\Sigma,\qquad
 u(0,\cdot)=u_0\quad\text{in }\Omega.
\end{equation}
For a stationary problem, $Q$ and $\Sigma$ are
replaced by $\Omega$ and $\partial\Omega$, and the initial condition is
omitted. We assume that the problem has a unique solution in the stated
solution class.

A PINN represents an approximate solution by a neural network $v_\theta$,
where $\theta\in\Theta\subset\mathbb R^P$ denotes the vector of trainable parameters and $P$ their number.
The input is $(t,x)$, or $x$ in a stationary problem, and the output represents
the unknown field or fields $u$. The derivatives appearing in the PDE are
computed by automatic differentiation of the network. The network must
therefore have the regularity required by the differential and boundary
operators. To train the network, we form a loss from the PDE and data residuals \citep{RaissiEtAl2019,DeRyckMishra2024}.

For an admissible approximation $v$, the residuals are defined by
\begin{equation}\label{eq:pinn-residuals}
 R(v):=\mathcal F(v)-f,\qquad
 B(v):=\mathcal B(v)-g,\qquad
 I(v):=v(0,\cdot)-u_0.
\end{equation}
If all these residuals are zero, then $v$
satisfies \eqref{eq:pinn-problem}, so uniqueness of the solution implies
$v=u$. The standard PINN loss therefore penalizes these residuals through a weighted sum of
squared $L^2$ norms given by
\begin{equation}\label{eq:standard-general}
 \Lstd(v):=\norm{R(v)}_{L^2(Q)}^2
       +\lambda_B\norm{B(v)}_{L^2(\Sigma)}^2
       +\lambda_I\norm{I(v)}_{L^2(\Omega)}^2,
 \qquad \lambda_B,\lambda_I>0.
\end{equation}
Here $\lambda_B$ and $\lambda_I$ denote the weights of the boundary and initial condition terms, respectively. Both weights are independent of $\theta$. A condition imposed exactly by the neural representation has zero
residual and its penalty can be omitted. Although implementations
approximate these integrals by numerical quadrature, all losses in our
analysis are defined by exact integrals.

To obtain an approximation $v_\theta$ of $u$, we seek parameters $\theta$
that approximately solve
\begin{equation}\label{eq:pinn-minimization}
 \inf_{\theta\in\Theta}\Lstd(v_\theta).
\end{equation}
Various optimization methods, including gradient descent,
are used to solve \eqref{eq:pinn-minimization}
\citep{RathoreEtAl2024,MullerZeinhofer2023}.

We also refer to \citep{RaissiEtAl2019,KarniadakisEtAl2021} for background
on PINNs and physics-informed machine learning, and to
\citep{DeRyckMishra2024} for their numerical analysis.

\subsection{Stability theory}\label{sec:stability-theory}
Stability describes how perturbations of a PDE and its data affect
the solution \citep{Evans1998,Brezis2011}. Limited network expressivity and incomplete optimization
can leave a positive loss after training. An approximation $v$ with
residuals \eqref{eq:pinn-residuals} satisfies the perturbed problem
\[
 \mathcal F(v)=f+R(v),\qquad
 \mathcal B(v)=g+B(v),\qquad
 v(0,\cdot)=u_0+I(v).
\]
Thus, a stability estimate for the PDE can convert residual bounds into
an error bound. This viewpoint is central to the error analysis of PINNs
\citep{DeRyckMishra2024,MishraMolinaro2023,DeRyckJagtapMishra2024,ZeinhoferEtAl2024}.

Let $\mathcal X$ be the normed space in which the error is measured, and
let $\mathcal A$ be the admissible class of functions satisfying the
regularity, a priori bounds and other conditions required by the
analysis. The exact solution satisfies $u\in\mathcal A$.
Writing $\E(v)$ for the error in $\mathcal X$, we consider stability estimates of the form
\begin{equation}\label{eq:standard-stability}
 (\E(v))^2:=\norm{v-u}_{\mathcal X}^2
 \leq C_{\mathrm{stab}}\Lstd(v),\qquad v\in\mathcal A,
\end{equation}
where the constant $C_{\mathrm{stab}}$ is independent of the individual
approximation $v$. If
$v_\theta\in\mathcal A$, the stability estimate~\eqref{eq:standard-stability} applies to the neural network approximation,
regardless of the algorithm that produced its parameters. In particular,
vanishing loss then implies convergence in $\mathcal X$.
The choice of $\mathcal X$ and the validity of \eqref{eq:standard-stability} must be
established for the particular PDE and residual norms.

\section{Error estimates with and without integral constraints}\label{sec:theory}
In this section, we use stability theory to explain when and why integral constraints
can improve the accuracy of PINN approximations. We consider a nonlinear
elliptic equation with Neumann boundary conditions, the Helmholtz and
stationary Burgers equations with Dirichlet boundary conditions, and a
time-dependent reaction--diffusion equation. In these four cases, the
standard loss can provide limited control of constant shifts of the
solution, oscillatory modes, shifts in the position of a viscous shock,
and a time-dependent spatial mean, respectively. Since such an error can
remain large even when the standard loss is small, the stability constant
in \eqref{eq:standard-stability} necessarily becomes large, as the factor
$\varepsilon^{-2}$ in \eqref{eq:intro-standard} illustrates. In these
regimes, reducing the standard loss alone cannot guarantee a target
accuracy. For each PDE,
we first quantify how the standard loss bound depends on the parameter
or final time and then show how adding an integral constraint derived
from the PDE and the given data improves error control.
For a spatial domain $\Omega$, we write $\norm{\cdot}_2:=\norm{\cdot}_{L^2(\Omega)}$,
$\norm{\cdot}_{2,\partial\Omega}:=\norm{\cdot}_{L^2(\partial\Omega)}$ and
$\norm{\cdot}_{H^1}:=\norm{\cdot}_{H^1(\Omega)}$.

\subsection{A nonlinear elliptic equation with Neumann boundary conditions}\label{sec:neumann}
We consider the following nonlinear elliptic PDE on a
bounded connected Lipschitz domain $\Omega\subset\R^d$ for $\eps>0$
\begin{equation}\label{eq:N}
 \begin{cases}
 -\Delta u+\eps g(u)=f,&x\in\Omega,\\
 \partial_nu=0,&x\in\partial\Omega,
 \end{cases}
\end{equation}
where $f\in L^2(\Omega)$ and $g\in C^1(I)$ is defined on an interval
$I\subset\R$ and satisfies
\begin{equation}\label{eq:g}
 0<\alpha\le g'(s)\le\beta<\infty\qquad(s\in I),
\end{equation}
with $\alpha$ and $\beta$ independent of $\eps$.
We assume that the Neumann problem~\eqref{eq:N} has a solution $u\in H^2(\Omega)$ taking values
in $I$.
For further details on the mathematical analysis of related PDEs,
we refer to \citep{Evans1998,GilbargTrudinger2001,Brezis2011}.

We consider approximations $v\in H^2(\Omega)$ taking values in $I$.
Under \eqref{eq:g}, $g(u)$ and $g(v)$ are in $L^2(\Omega)$.
Following the standard PINN formulation, we define the residuals and standard loss by
\begin{equation}\label{eq:Nres}
 r:=-\Delta v+\eps g(v)-f,\qquad b:=\partial_nv,\qquad
 \mathcal L_{\mathrm{std},\eps}(v):=\norm r_{L^2(\Omega)}^2+
 \norm b_{L^2(\partial\Omega)}^2.
\end{equation}
For a constant $c$ such that $u+c$ takes values in $I$, substituting
$v=u+c$ in \eqref{eq:Nres} gives $b=0$ and
$r=\eps[g(u+c)-g(u)]$. Only the reaction term detects the shift, and
the upper bound in~\eqref{eq:g} gives $|r|\le\eps\beta|c|$. Thus a fixed nonzero shift can
produce a small standard loss when $\eps$ is small.

The following
theorem quantifies the dependence of the $H^1(\Omega)$ error bound
on $\eps$.

\begin{theorem}[Standard loss estimate]\label{prop:Nstandard}
For every $v\in H^2(\Omega)$ taking values in $I$, there is a constant $C$,
depending only on $\Omega$ and $\alpha$, such that for all $0<\eps\le1$,
\begin{equation}\label{eq:Nstandard}
 \norm{v-u}_{H^1}^2\le C\eps^{-2}
 \mathcal L_{\mathrm{std},\eps}(v).
\end{equation}
In fact, if there is a constant $c\ne0$, independent of $\eps$, such that
$u+c$ takes values in $I$ for all sufficiently small $\eps>0$, then
$\tilde v:=u+c$ satisfies
\begin{equation}\label{eq:Nexample}
 \norm{\tilde v-u}_{H^1}=|c||\Omega|^{1/2}>0,\qquad
 \mathcal L_{\mathrm{std},\eps}(\tilde v)\le\beta^2c^2|\Omega|\eps^2.
\end{equation}
Consequently, the $\eps^{-2}$ dependence in \eqref{eq:Nstandard} is sharp.
\end{theorem}
\begin{proof}
$\bullet$ \emph{On the estimate~\eqref{eq:Nstandard}.}
Let $e:=v-u$ and define
\begin{equation}\label{eq:Ncoefficient}
 a(x):=\int_0^1g'\bigl(u(x)+t e(x)\bigr)\dd t.
\end{equation}
Then $g(v)-g(u)=ae$ by the fundamental theorem of calculus, and $\alpha\le a\le\beta$ by~\eqref{eq:g}.
Subtracting the PDE and boundary condition in~\eqref{eq:N}
from the corresponding equations for $v$ in~\eqref{eq:Nres} and using
the identity $g(v)-g(u)=ae$ gives
\begin{equation}\label{eq:Nerror}
 \begin{cases}
 -\Delta e+\eps ae=r,&x\in\Omega,\\
 \partial_n e=b,&x\in\partial\Omega.
 \end{cases}
\end{equation}
Multiplying~\eqref{eq:Nerror} by $e$ and integrating by parts gives
\[
 \norm{\nabla e}_2^2+\eps\int_\Omega ae^2
 =\int_\Omega re\dd x+\int_{\partial\Omega}be\dd s
 \le C(\norm r_2+\norm b_{2,\partial\Omega})\norm e_{H^1}.
\]
Since $0<\eps\le1$ and $a\ge\alpha$, we have
\[
 \begin{aligned}
 \norm{\nabla e}_2^2+\eps\int_\Omega ae^2
 &\ge\eps\norm{\nabla e}_2^2+\eps\alpha\norm e_2^2\\
 &\ge\eps\min\{1,\alpha\}\norm e_{H^1}^2.
 \end{aligned}
\]
If $\norm e_{H^1}=0$, the inequality~\eqref{eq:Nstandard} is immediate. Otherwise, dividing the
preceding energy inequality by
$\eps\min\{1,\alpha\}\norm e_{H^1}$ gives
\[
 \norm e_{H^1}\le
 \frac{C}{\eps\min\{1,\alpha\}}
 (\norm r_2+\norm b_{2,\partial\Omega}).
\]
It follows that
\[
 \begin{aligned}
 \norm e_{H^1}^2
 &\le \frac{2C^2}{\eps^2\min\{1,\alpha\}^2}
 \bigl(\norm r_2^2+\norm b_{2,\partial\Omega}^2\bigr)\\
 &\le C\eps^{-2}\mathcal L_{\mathrm{std},\eps}(v),
 \end{aligned}
\]
where we used the loss definition~\eqref{eq:Nres} and
absorbed the fixed factors into the constant $C$ of Theorem~\ref{prop:Nstandard}.

\smallskip
\noindent$\bullet$ \emph{On the existence of $\tilde v$.}
Under the additional assumption, the shift
$\tilde v=u+c$ has zero boundary residual and satisfies
$|g(u+c)-g(u)|\le\beta|c|$, giving \eqref{eq:Nexample}.
Its squared error divided by the standard loss is therefore at least
$\beta^{-2}\eps^{-2}$, proving the asserted sharpness in
\eqref{eq:Nstandard}.
\end{proof}
For any $\delta>0$, the estimate~\eqref{eq:Nstandard} guarantees
$\norm{v-u}_{H^1}\le\delta$ whenever
$\mathcal L_{\mathrm{std},\eps}(v)\le\eps^2\delta^2/C$.
This sufficient condition sets an extremely small loss threshold when
$\eps$ is small. In \eqref{eq:Nexample}, the approximation $\tilde v=u+c$ has exactly
the same shape as $u$, although its $H^1$ error stays bounded away from
zero. Thus, at small $\eps$, a network with a small standard loss may
reproduce the shape of $u$ with a nearly constant offset throughout
$\Omega$.

On the other hand, integrating \eqref{eq:N} over $\Omega$ and using the Neumann
boundary condition gives
\begin{equation}\label{eq:neumann-identity}
 \int_\Omega g(u)\dd x=\eps^{-1}\int_\Omega f\dd x.
\end{equation}
Using the integral $\int_\Omega g(u)\dd x$ obtained in this way, we define
the following augmented loss so that the approximation matches this value:
\begin{equation}\label{eq:neumann-aug}
 \mathcal L_{\mathrm{aug},\eps}(v)
 :=\mathcal L_{\mathrm{std},\eps}(v)+|D(v)|^2,
\end{equation}
where
\[
 D(v):=\int_\Omega g(v)\dd x-\int_\Omega g(u)\dd x.
\]
For this loss, we obtain the following estimate.

\begingroup
\renewcommand{\pinnTheoremName}{Theorem}
\begin{theorem}[Uniform estimate with an integral constraint]\label{thm:N}
For every $v\in H^2(\Omega)$ taking values in $I$, there is a constant $C$,
depending only on $\Omega$, $\alpha$ and $\beta$, such that for all $\eps>0$,
\begin{equation}\label{eq:Nuniform}
 \norm{v-u}_{H^1(\Omega)}^2
 \le C\mathcal L_{\mathrm{aug},\eps}(v).
\end{equation}
\end{theorem}
\endgroup
\begin{proof}
With $e=v-u$ and $a$ as in \eqref{eq:Ncoefficient}, we have
$D(v)=\int_\Omega ae$.
We define $c_a$ and $w$ by
\[
 c_a:=\frac{D(v)}{\int_\Omega a\dd x},\qquad w:=e-c_a.
\]
We have $\int_\Omega aw=0$ and $|c_a|\le |D(v)|/(\alpha|\Omega|)$.
Let $\overline w:=|\Omega|^{-1}\int_\Omega w\dd x$ denote the spatial mean of $w$. Since $\int_\Omega aw=0$,
\[
 \overline w\int_\Omega a=-\int_\Omega a(w-\overline w).
\]
The Poincar\'e inequality then implies
\begin{equation}\label{eq:weightedP}
 \norm w_2\le(1+\beta/\alpha)\norm{w-\overline w}_2
 \le C(1+\beta/\alpha)\norm{\nabla w}_2.
\end{equation}
Multiplying~\eqref{eq:Nerror} by $w$ and integrating by parts,
$\nabla e=\nabla w$ and $\int_\Omega aw\dd x=0$ imply that the cross
term $\eps c_a\int_\Omega aw\dd x$ vanishes, giving
\begin{equation}\label{eq:Nenergy}
 \norm{\nabla w}_2^2+\eps\int_\Omega aw^2
 =\int_\Omega rw\dd x+\int_{\partial\Omega}bw\dd s.
\end{equation} By \eqref{eq:weightedP} and the trace inequality,
the right-hand side of \eqref{eq:Nenergy} is at most
$C(\norm r_2+\norm b_{2,\partial\Omega})\norm{\nabla w}_2$.
Dropping the nonnegative reaction term gives
\[
 \norm w_{H^1}\le C(\norm r_2+\norm b_{2,\partial\Omega}).
\]
Since $c_a$ is constant in space,
\[
 \norm{c_a}_{H^1}=|c_a||\Omega|^{1/2}
 \le\frac{|D(v)|}{\alpha|\Omega|^{1/2}}.
\]
Recalling that $w=e-c_a$, we obtain
\[
 \begin{aligned}
 \norm e_{H^1}
 &\le\norm w_{H^1}+\norm{c_a}_{H^1}\\
 &\le C\bigl(\norm r_2+\norm b_{2,\partial\Omega}+|D(v)|\bigr).
 \end{aligned}
\]
It follows that
\[
 \norm e_{H^1}^2
 \le C\bigl(\norm r_2^2+\norm b_{2,\partial\Omega}^2+|D(v)|^2\bigr)
 =C\mathcal L_{\mathrm{aug},\eps}(v).
\]
The constant $C$ depends only on $\Omega$, $\alpha$ and $\beta$, and is
independent of $\eps$. This proves \eqref{eq:Nuniform}.
\end{proof}
We take $C$ to be the larger of the constants in Theorems~\ref{prop:Nstandard} and~\ref{thm:N}. For $0<\eps\le1$, the estimates~\eqref{eq:Nstandard}
and~\eqref{eq:Nuniform} give the same error bound
$\norm{v-u}_{H^1}\le\delta$ under the respective sufficient conditions
\[
 \mathcal L_{\mathrm{std},\eps}(v)
 \le \frac{\eps^2\delta^2}{C},\qquad
 \mathcal L_{\mathrm{aug},\eps}(v)
 \le \frac{\delta^2}{C}.
\]
This constant is independent of $\eps$, so the threshold
for the augmented loss is independent of $\eps$ while that for the standard
loss decreases in proportion to $\eps^2$.
When the two losses are reduced to comparable values, the integral
constraint controls constant shifts uniformly in $\eps$ and is therefore
expected to improve accuracy for small $\eps$.
When $\eps$ is near one, however, the factor $\eps^{-2}$ in
\eqref{eq:Nstandard} is close to one, so adding the integral constraint may give only a small gain in accuracy.

\begin{remark}[Relation to the pure Neumann problem]\label{rem:pure-neumann}
For $\eps=0$, equation \eqref{eq:N} is solvable only when
$\int_\Omega f\dd x=0$, and it then determines $u$ only up to an additive
constant, so the mean of the solution has to be given for the problem
to be uniquely solvable. With a penalty on the mean of the approximation,
the $H^1$ error is controlled by the excess of the Deep Ritz energy over
its minimum \citep{LuLuWang2021}, which is also presented in
\citep[Example 5.11]{DeRyckMishra2024}, and by the strong residual loss of
mixed form in \citep{BaiEtAl2025}. In that setting, the mean condition is part
of the problem itself. Here $\eps>0$, so the Neumann problem~\eqref{eq:N} is uniquely solvable
without any additional condition, and the identity~\eqref{eq:neumann-identity} is a consequence
of the PDE and the boundary condition. Theorem~\ref{prop:Nstandard}
shows that the standard loss controls constant shifts only weakly when
$\eps$ is small, and Theorem~\ref{thm:N} shows that adding this redundant
information through the integral constraint in \eqref{eq:neumann-aug}
restores a bound uniform in $\eps$. The reason for including the integral
constraint therefore differs from the reason in the pure Neumann problem.
\end{remark}

\begin{remark}[Pseudo-Two-Dimensional model in Lithium-Ion Batteries]\label{rem:battery}
At each fixed time, the positive-electrode solid-phase potential $\psi_p$
in the P2D model satisfies \citep{LeeEtAl2025}
\begin{equation}\label{eq:battery-potential}
 \begin{cases}
 \sigma_p^{\rm eff}\partial_{xx}\psi_p=a_pFj_p,&0<x<L_p,\\
 -\sigma_p^{\rm eff}\partial_x\psi_p=I_{\rm app},&x=0,\\
 -\sigma_p^{\rm eff}\partial_x\psi_p=0,&x=L_p.
 \end{cases}
\end{equation}
Here $L_p$ denotes the electrode thickness, $\sigma_p^{\rm eff}$ the effective electronic conductivity, $a_p$ the active surface area per unit electrode volume, and $I_{\rm app}$ the applied current density. The Butler--Volmer ionic flux $j_p$ is proportional to
$\sinh(F(\psi_p-\phi_p-U_p)/(2RT))$, with a concentration-dependent
prefactor. In this expression, $\phi_p$ and $U_p$ denote the liquid-phase and open circuit potentials, $T$ denotes temperature, and $F$ and $R$ denote the Faraday and gas constants, respectively. Integrating the differential equation in \eqref{eq:battery-potential}
over $(0,L_p)$ and using its boundary conditions gives the normalized
secondary conservation law
\begin{equation}\label{eq:battery-integral}
 {a_pF\over I_{\rm app}}\int_0^{L_p}j_p\,dx=1
 \qquad(I_{\rm app}\ne0).
\end{equation}

In \citep{LeeEtAl2025}, minimizing the PINN loss did not provide accurate
predictions for the P2D model, even when a small loss was achieved.
Adding integral constraints such as \eqref{eq:battery-integral} to the loss
improved the prediction accuracy. With all other state variables exact,
the potential $\psi_p$ in~\eqref{eq:battery-potential} satisfies a nonlinear elliptic equation with Neumann
boundary conditions in which the reaction term depends on $x$ through those
variables and the flux is given at $x=0$. If the reaction derivative
satisfies bounds of the form \eqref{eq:g} uniformly in $x$, the arguments
of this subsection extend to this setting. Thus, the improvement in
prediction accuracy from adding integral constraints such as
\eqref{eq:battery-integral} can be explained in the same way.
\end{remark}

\subsection{Helmholtz equation with Dirichlet boundary conditions}\label{sec:helmholtz}
In the preceding Neumann example, constant shifts leave the Laplacian and
the Neumann data unchanged, so they enter the residual only through the
reaction term. The given Dirichlet values rule out such shifts, but a
small residual can coexist with a large solution error through other
mechanisms. The following Helmholtz example illustrates one such mechanism.

For $0<\eps\le1$, we consider the following Helmholtz equation
\begin{equation}\label{eq:HelmPDE}
 \begin{cases}
 -u''-(4\pi^2-\eps)u=f,&0<x<1,\\
 u(0)=u(1)=0.
 \end{cases}
\end{equation}
For every $f\in L^2(0,1)$, the Dirichlet
problem~\eqref{eq:HelmPDE} has a unique solution
$u\in H^2(0,1)\cap H^1_0(0,1)$. We note that if we assume $\eps=0$ and $u_0$ is a solution,
then so is $u_0+c\sin(2\pi x)$ for every $c\in\R$, since
$\sin(2\pi x)$ vanishes at both endpoints and satisfies
$-\partial_x^2\sin(2\pi x)=4\pi^2\sin(2\pi x)$.
Further details on existence and uniqueness can be found in \citep[Chapter~6]{Evans1998}.
Related applications in wave propagation and acoustic scattering can be
found in \citep{Ihlenburg1998}.

This nonuniqueness at $\eps=0$ helps explain what happens at a small
positive $\eps$. Adding $c\sin(2\pi x)$ to the exact solution then changes
the PDE residual by only $\eps c\sin(2\pi x)$, so a small residual can
coexist with a large solution error. We now ask how large the $L^2$ error can be for a given residual norm.
For $v\in H^2(0,1)\cap H^1_0(0,1)$, we define the residual and the standard loss by
\begin{equation}\label{eq:HelmLoss}
 r:=-v''-(4\pi^2-\eps)v-f,\qquad
 \mathcal L_{\mathrm{std},\eps}(v):=\norm r_{L^2(0,1)}^2.
\end{equation}
The following theorem gives the largest possible $L^2$ error for a given residual norm.

\begin{theorem}[Standard loss estimate]\label{prop:HelmStandard}
For every $v\in H^2(0,1)$ with $v(0)=v(1)=0$ and every
$0<\eps\le1$,
\begin{equation}\label{eq:HelmStandard}
 \norm{v-u}_2^2\le\eps^{-2}\mathcal L_{\mathrm{std},\eps}(v).
\end{equation}
Here, $\eps^{-2}$ is the exact best constant. In fact, for each $0<\eps\le1$,
there is $\tilde v\in H^2(0,1)\cap H^1_0(0,1)$ such that
\begin{equation}\label{eq:HelmExample}
 \norm{\tilde v-u}_2=1,\qquad
 \mathcal L_{\mathrm{std},\eps}(\tilde v)=\eps^2.
\end{equation}
\end{theorem}
\begin{proof}
$\bullet$ \emph{On the estimate~\eqref{eq:HelmStandard}.}
Let us define $\phi_n(x):=\sqrt2\sin(n\pi x)$,
$n\ge1$. These Dirichlet eigenfunctions
form an orthonormal basis of $L^2(0,1)$ and satisfy
$-\phi_n''=(n\pi)^2\phi_n$.

Let $e:=v-u$ and define its sine coefficients by
\[
 e_n:=\int_0^1 e(x)\phi_n(x)\dd x,\qquad n\ge1.
\]
Since $u$ satisfies
\eqref{eq:HelmPDE}, we have $r=-e''-(4\pi^2-\eps)e$.
Integrating by parts twice and using
$e(0)=e(1)=\phi_n(0)=\phi_n(1)=0$, we obtain
\[
 \begin{aligned}
 \int_0^1(-e'')\phi_n\dd x
 &=\bigl[-e'\,\phi_n+e\,\phi_n'\bigr]_0^1
     -\int_0^1e\,\phi_n''\dd x\\
 &=(n\pi)^2e_n.
 \end{aligned}
\]
Let $r_n:=\int_0^1r\phi_n\dd x$ denote the $n$th sine coefficient of $r$. The preceding identity gives
\[
 r_n=\bigl((n^2-4)\pi^2+\eps\bigr)e_n.
\]
Applying Parseval's identity to $r$ now gives
\[
 \norm r_2^2=\sum_{n\ge1}|r_n|^2
 =\sum_{n\ge1}\bigl((n^2-4)\pi^2+\eps\bigr)^2|e_n|^2.
\]
The coefficient inside the square equals $\eps$ for $n=2$. For $n\ne2$,
$|(n^2-4)\pi^2+\eps|\ge3\pi^2-1>1\ge\eps$, since $0<\eps\le1$.
Applying Parseval's identity also to $e$, we therefore obtain
\[
 \norm r_2^2\ge\eps^2\sum_{n\ge1}|e_n|^2
 =\eps^2\norm{v-u}_2^2,
\]
which gives \eqref{eq:HelmStandard} after division by $\eps^2$.

\smallskip
\noindent$\bullet$ \emph{On the existence of $\tilde v$.}
For $\tilde v:=u+c\phi_2$, we have
$\norm{\tilde v-u}_2^2=c^2$ and
$\mathcal L_{\mathrm{std},\eps}(\tilde v)=\eps^2c^2$.
Thus equality is attained, and taking $c=1$ gives \eqref{eq:HelmExample}.
\end{proof}
For $v:=u+\sin(2\pi x)$, the standard loss is
$\eps^2/2$ even though $\norm{v-u}_2=1/\sqrt{2}$. The residual
controls all other sine modes uniformly in $\eps$, so any appreciable error
at very small standard loss must be close in $L^2$ to a multiple of
$\sin(2\pi x)$.

To measure this error component, we multiply \eqref{eq:HelmPDE} by $\phi_2$ and
integrate by parts twice. The Dirichlet boundary values give
\[
 \int_0^1 u\phi_2\dd x=\eps^{-1}\int_0^1 f\phi_2\dd x.
\]
As in Section~\ref{sec:neumann}, using the integral $\int_0^1 u\phi_2\dd x$ obtained in
this way, we define the following augmented loss so that the approximation
matches this value:
\begin{equation}\label{eq:helmholtz-aug}
 \mathcal L_{\mathrm{aug},\eps}(v)
 :=\mathcal L_{\mathrm{std},\eps}(v)+|D(v)|^2,
\end{equation}
where
\[
 D(v):=\int_0^1 v\phi_2\dd x-\int_0^1 u\phi_2\dd x.
\]
For this loss, we obtain the following estimate.

\begin{theorem}[Uniform estimate with a weighted integral constraint]\label{prop:HelmUniform}
For every $v\in H^2(0,1)$ with $v(0)=v(1)=0$ and every
$0<\eps\le1$,
\begin{equation}\label{eq:HelmUniform}
 \norm{v-u}_2^2
 \le (3\pi^2-1)^{-2}\norm r_2^2+|D(v)|^2
 \le\mathcal L_{\mathrm{aug},\eps}(v).
\end{equation}
\end{theorem}
\begin{proof}
Recalling $D(v)=\int_0^1(v-u)\phi_2\dd x$, we have $D(v)=e_2$.
The bound $|(n^2-4)\pi^2+\eps|\ge3\pi^2-1$ for $n\ne2$ gives
\[
 \begin{aligned}
 \sum_{n\ne2}|e_n|^2
 &\le(3\pi^2-1)^{-2}
 \sum_{n\ne2}\bigl((n^2-4)\pi^2+\eps\bigr)^2|e_n|^2\\
 &\le(3\pi^2-1)^{-2}\norm r_2^2.
 \end{aligned}
\]
Since $(3\pi^2-1)^{-2}<1$, Parseval's identity now yields
\[
 \begin{aligned}
 \norm{v-u}_2^2
 &=|e_2|^2+\sum_{n\ne2}|e_n|^2\\
 &\le|D(v)|^2+(3\pi^2-1)^{-2}\norm r_2^2\\
 &\le|D(v)|^2+\norm r_2^2
 =\mathcal L_{\mathrm{aug},\eps}(v).
 \end{aligned}
\]
This proves \eqref{eq:HelmUniform}.
\end{proof}
The integral in~\eqref{eq:helmholtz-aug} measures the $\phi_2$ component of
the error directly. At the same loss value $\ell$, the standard
estimate~\eqref{eq:HelmStandard} and the augmented estimate~\eqref{eq:HelmUniform}
give upper bounds of $\ell/\eps^2$ and $\ell$, respectively, for the squared
$L^2$ error. The bounds agree at $\eps=1$ and suggest a greater
gain in accuracy at small $\eps$ when training attains comparable losses.

\begin{remark}
For simplicity, we focus on the one-dimensional Helmholtz problem. On a smooth bounded domain in
any dimension, the same argument applies with $4\pi^2$ replaced by a
Dirichlet eigenvalue of $-\Delta$ and $\eps>0$ sufficiently small.
If the eigenvalue has multiplicity greater than one, we use one integral
constraint for each eigenfunction in an orthonormal basis of its eigenspace.
\end{remark}

\subsection[The stationary viscous Burgers equation]{The stationary viscous Burgers equation}\label{sec:burgers}
We next consider the stationary viscous Burgers equation with the same
convective flux at both boundary endpoints. We note that its exact solution
has a very narrow transition layer at small viscosity.
In Section~\ref{sec:neumann}, a shift of the solution graph along the
$y$-axis can produce a large error at a very small standard loss.
In the present example, a large error at a very small standard loss can
arise from a shift of the transition layer along the $x$-axis.
This example illustrates how an integral constraint can improve the
accuracy of PINNs by controlling this shift of the transition layer.

Specifically, we study the following boundary value problem
\begin{equation}\label{eq:Bproblem}
 \begin{cases}
 -\eps u''+u\,u'=0, & x\in I:=(-1,1),\\
 u(-1)=1,\qquad u(1)=-1,
 \end{cases}
\end{equation}
where $0<\eps\le1$. The asymptotic behavior for small viscosity and the existence of stationary
solutions are studied in \citep{GarbeyKaper1999,FogelklouEtAl2008}.
For more details on the viscous Burgers equation, we refer to
\citep{KreissKreiss1986,LaforgueOMalley1994,ReynaWard1995}.

The exact solution can be written explicitly as
\begin{equation}\label{eq:Bsolution}
 u(x)=-A\tanh\frac{Ax}{2\eps},
\end{equation}
where $A\in(1,2)$ is the unique positive root of
\[
 A\tanh\frac{A}{2\eps}=1.
\]

For every $v\in H^2(I)$ with $v(-1)=1$ and $v(1)=-1$,
we define the residual and standard loss by
\begin{equation}\label{eq:Bres}
 r:=-\eps v''+v\,v',\qquad
 \Lstd(v):=\norm r_2^2.
\end{equation}
The following theorem shows that the optimal constant in the bound
for the squared $L^2(I)$ error in terms of the standard loss \eqref{eq:Bres}
is of order $\eps e^{2/\eps}$ as $\eps\to0$.

\begin{theorem}[Standard loss estimate]\label{prop:Bstandard}
For every $v\in H^2(I)$ with $v(-1)=1$ and $v(1)=-1$,
there is a constant $C$, independent of $v$ and $\eps$, such
that for every $0<\eps\le1$,
\begin{equation}\label{eq:Bstandard}
 \norm{v-u}_2^2\le C\eps e^{2/\eps}\Lstd(v).
\end{equation}
Here, the $\eps e^{2/\eps}$ dependence is sharp as $\eps\to0$.
For each sufficiently small $\eps>0$, there is
$\tilde v\in H^2(I)$ with $\tilde v(-1)=1$ and $\tilde v(1)=-1$ such that
\begin{equation}\label{eq:Bexample}
 \norm{\tilde v-u}_2\ge\frac12,\qquad
 \Lstd(\tilde v)\le C e^{-3/(2\eps)}.
\end{equation}
\end{theorem}
\begin{proof}
$\bullet$ \emph{On the upper bound~\eqref{eq:Bstandard}.}
If $v\equiv u$, the estimate \eqref{eq:Bstandard} is immediate.
We therefore assume $v\not\equiv u$ and consider two cases according
to the size of $\norm r_2$. For small $\norm r_2$, we decompose $u-v$
into two terms and estimate them separately. For large $\norm r_2$,
we estimate $\norm{u-v}_2$ directly.

Let us begin with the case $\norm r_2\le 1/(4\sqrt2)$. Since
$H^2(I)\hookrightarrow C^1(\bar I)$ in one dimension, we use the
$C^1(\bar I)$ representative of $v$. The boundary values $v(-1)=1$
and $v(1)=-1$ then imply that $v$ has a zero in $I$.
We fix one such zero $s$ and define the translate $u_s$ of
the exact profile \eqref{eq:Bsolution} by
\begin{equation}\label{eq:Bprofile}
 u_s(x):=-A\tanh\frac{A(x-s)}{2\eps},\qquad x\in\R.
\end{equation}
With this $u_s$, we decompose the error as
\begin{equation}\label{eq:Bsplit}
 u-v=(u-u_s)+(u_s-v).
\end{equation}
We first estimate the second term $u_s-v$ in \eqref{eq:Bsplit}.
To obtain upper and lower bounds for $v$, we derive differential inequalities for $v$.
\begin{figure}[!t]
\centering
\includegraphics[width=0.94\linewidth]{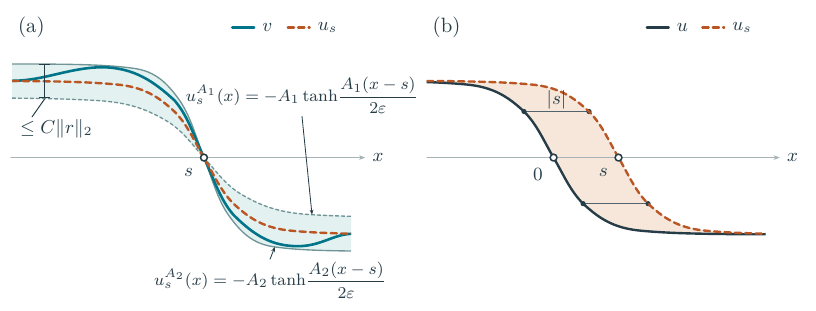}
\caption{Sketches of the estimates for (a) $v-u_s$ and (b) $u_s-u$.}
\label{fig:burgers-decomposition}
\end{figure}

Since $v\not\equiv u$, we choose $x_0\in I$ with $v(x_0)\ne u(x_0)$.
By continuity, the zeros of $v-u$ in $[-1,x_0]$ form a compact set.
Since $v(-1)=u(-1)$, there is at least one zero in $[-1,x_0]$, so
a largest zero exists. Likewise, the zeros of $v-u$ in $[x_0,1]$
form a compact set containing $1$, so a smallest zero exists. We write
\[
 x_-:=\max\{x\in[-1,x_0]:v(x)=u(x)\},\qquad
 x_+:=\min\{x\in[x_0,1]:v(x)=u(x)\}.
\]
Since $v(x_0)\ne u(x_0)$, we have $x_-<x_0<x_+$. By the choice of $x_-$,
there is no zero in $(x_-,x_0)$, and by the choice of $x_+$, there is no
zero in $(x_0,x_+)$. Hence $v-u$ has the same sign throughout $(x_-,x_+)$.
If $v>u$ on $(x_-,x_+)$, the equalities $v(x_-)=u(x_-)$ and $v(x_+)=u(x_+)$ give
$(v-u)'(x_-)\ge0$ and $(v-u)'(x_+)\le0$. If $v<u$ on $(x_-,x_+)$, then
$(v-u)'(x_-)\le0$ and $(v-u)'(x_+)\ge0$.

With these facts, we evaluate $-\eps v'+v^2/2$ at $x_-$ and $x_+$.
Differentiating the exact solution \eqref{eq:Bsolution}, we obtain
\begin{equation}\label{eq:BexactODE}
 -\eps u'(x)+\frac{(u(x))^2}{2}
 =\frac{A^2}{2}\left(\operatorname{sech}^2\frac{Ax}{2\eps}
                     +\tanh^2\frac{Ax}{2\eps}\right)
 =\frac{A^2}{2}.
\end{equation}
Here, for $x=x_-$ or $x=x_+$, we substitute \eqref{eq:BexactODE} for $A^2/2$
and then use $v(x)=u(x)$ to obtain
\begin{equation}\label{eq:BdifferenceODE}
 \begin{aligned}
 -\eps v'(x)+\frac{(v(x))^2}{2}-\frac{A^2}{2}
 &= -\eps v'(x)+\frac{(v(x))^2}{2}
    -\left(-\eps u'(x)+\frac{(u(x))^2}{2}\right)\\
 &= -\eps(v-u)'(x)+\frac{(v(x))^2-(u(x))^2}{2}\\
 &= -\eps(v-u)'(x).
 \end{aligned}
\end{equation}
If $v>u$ on $(x_-,x_+)$, the signs of $(v-u)'(x_-)$ and $(v-u)'(x_+)$
in \eqref{eq:BdifferenceODE} give
\[
 -\eps v'(x_-)+\frac{(v(x_-))^2}{2}
 \le\frac{A^2}{2}\le -\eps v'(x_+)+\frac{(v(x_+))^2}{2}.
\]
If $v<u$ on $(x_-,x_+)$, both inequalities are reversed. In either case,
the continuity of $-\eps v'+v^2/2$ gives a point $y\in[x_-,x_+]$ such that
\[
 -\eps v'(y)+\frac{(v(y))^2}{2}=\frac{A^2}{2}.
\]

We note from \eqref{eq:Bres} that $(-\eps v'+v^2/2)'=r$.
Integrating this identity from $y$ to $x\in\bar I$, we obtain
\begin{equation}\label{eq:Bfluxvariation}
 -\eps v'(x)+\frac{(v(x))^2}{2}-\frac{A^2}{2}
 =\int_y^x r(z)\dd z.
\end{equation}
Since $|\int_y^x r(z)\dd z|\le |x-y|^{1/2}\norm r_2\le\sqrt2\,\norm r_2$
by the Cauchy--Schwarz inequality, we obtain from \eqref{eq:Bfluxvariation}
\begin{equation}\label{eq:Bflux}
 \frac{A^2}{2}-\sqrt2\,\norm r_2
 \le -\eps v'(x)+\frac{(v(x))^2}{2}
 \le \frac{A^2}{2}+\sqrt2\,\norm r_2
 \qquad (x\in\bar I).
\end{equation}
We define
\[
 A_1:=\sqrt{A^2-2\sqrt2\,\norm r_2},\qquad
 A_2:=\sqrt{A^2+2\sqrt2\,\norm r_2}.
\]
Since $A>1$ and $\norm r_2\le1/(4\sqrt2)$, we have $A_1^2\ge1/2$.
Solving the inequalities
\eqref{eq:Bflux} for $\eps v'$ gives
\begin{equation}\label{eq:Bdiffineq}
 \frac{v^2-A_2^2}{2}\le\eps v'\le\frac{v^2-A_1^2}{2}.
\end{equation}
To obtain bounds for $v(x)$ from \eqref{eq:Bdiffineq}, we use
the formula for the translated profile $u_s$.
Translating the identity \eqref{eq:BexactODE} by $s$, we also have
$-\eps u_s'+(u_s)^2/2=A^2/2$. For $a>0$, we replace $A$ by $a$
in \eqref{eq:Bprofile} and define
\begin{equation}\label{eq:BprofileODE}
 u_s^a(x):=-a\tanh\frac{a(x-s)}{2\eps},\qquad x\in\R.
\end{equation}
It follows directly from \eqref{eq:BprofileODE} that
\[
 \eps(u_s^a)'=\frac{(u_s^a)^2-a^2}{2},\qquad u_s^a(s)=0.
\]
In particular, for $a=A_1$ and $a=A_2$, we obtain the equations
\[
 \eps(u_s^{A_1})'=\frac{(u_s^{A_1})^2-A_1^2}{2},\qquad
 \eps(u_s^{A_2})'=\frac{(u_s^{A_2})^2-A_2^2}{2}.
\]
Subtracting the equation for $u_s^{A_1}$ from the upper bound for
$\eps v'$ in \eqref{eq:Bdiffineq} and the equation for $u_s^{A_2}$
from the lower bound, we obtain
\begin{align}
 \eps(v-u_s^{A_1})'
 &\le\frac{v^2-(u_s^{A_1})^2}{2}
 =\frac{v+u_s^{A_1}}{2}(v-u_s^{A_1}),\label{eq:BdiffAone}\\
 \eps(v-u_s^{A_2})'
 &\ge\frac{v^2-(u_s^{A_2})^2}{2}
 =\frac{v+u_s^{A_2}}{2}(v-u_s^{A_2}).\label{eq:BdiffAtwo}
\end{align}
To show that $v(x)$ lies between $u_s^{A_1}(x)$ and $u_s^{A_2}(x)$,
we use the positive integrating factor
\[
 J_a(x):=\exp\!\left(-\frac{1}{2\eps}\int_s^x(v+u_s^a)(t)\dd t\right)>0,
 \qquad a\in\{A_1,A_2\}.
\]
For $s\le x\le1$, we multiply \eqref{eq:BdiffAone} by $J_{A_1}/\eps$
and \eqref{eq:BdiffAtwo} by $J_{A_2}/\eps$ and integrate from $s$ to $x$
to obtain
\begin{align}
 J_{A_1}(x)\bigl(v(x)-u_s^{A_1}(x)\bigr)
 &=\int_s^x [J_{A_1}(v-u_s^{A_1})]'(t)\dd t\le0,\label{eq:BproductAone}\\
 J_{A_2}(x)\bigl(v(x)-u_s^{A_2}(x)\bigr)
 &=\int_s^x [J_{A_2}(v-u_s^{A_2})]'(t)\dd t\ge0.\label{eq:BproductAtwo}
\end{align}
Here we used $v(s)=u_s^{A_1}(s)=u_s^{A_2}(s)=0$ at the lower endpoint
and $\eps J_a'=-(v+u_s^a)J_a/2$ to write the integrand as
\[
 [J_a(v-u_s^a)]'
 =\frac{J_a}{\eps}\left(\eps(v-u_s^a)'-\frac{v+u_s^a}{2}(v-u_s^a)\right).
\]
For $-1\le x\le s$, the inequalities in \eqref{eq:BproductAone}
and \eqref{eq:BproductAtwo} reverse since $\int_s^x=-\int_x^s$.
Dividing by the positive factors $J_{A_1}(x)$ and $J_{A_2}(x)$, we obtain
\begin{equation}\label{eq:Bcompare}
 \begin{aligned}
  u_s^{A_2}(x)&\le v(x)\le u_s^{A_1}(x) &&(s\le x\le1),\\
  u_s^{A_1}(x)&\le v(x)\le u_s^{A_2}(x) &&(-1\le x\le s).
 \end{aligned}
\end{equation}
By \eqref{eq:BprofileODE} and \eqref{eq:Bcompare}, we have
\[
 v(x)\ge u_s^{A_1}(x)>0\quad(-1\le x<s),\qquad
 v(x)\le u_s^{A_1}(x)<0\quad(s<x\le1),
\]
so $s$ is the only zero of $v$.

As shown in Figure~\ref{fig:burgers-decomposition}(a), $v(x)$ lies between
$u_s^{A_1}(x)$ and $u_s^{A_2}(x)$. Since $u_s^a(x)$ is monotone in $a$ for
each fixed $x$ by \eqref{eq:BprofileODE} and $A_1\le A\le A_2$, the band
also contains $u_s(x)=u_s^A(x)$. We can therefore estimate $|v(x)-u_s(x)|$
by the larger of the vertical distances from $u_s(x)$ to the two boundaries.
Indeed, we have
\[
 |v(x)-u_s(x)|\le
 \max_{i=1,2}|u_s^{A_i}(x)-u_s(x)|.
\]
We differentiate \eqref{eq:BprofileODE} with respect to $a$,
so that we have
\[
 \partial_a u_s^a(x)=-\tanh\!\left(\frac{a(x-s)}{2\eps}\right)
 -\frac{a(x-s)}{2\eps}\operatorname{sech}^2\!\left(\frac{a(x-s)}{2\eps}\right).
\]
For every $z\in\R$, we have $\cosh^2 z=1+\sinh^2 z\ge1+z^2$, so
$|z|\operatorname{sech}^2 z\le |z|/(1+z^2)\le 1/2$. Together with
$|\tanh z|\le1$, this gives $|\partial_a u_s^a(x)|\le 3/2\le2$.
For $i=1,2$, the mean value theorem gives
\[
 |u_s^{A_i}(x)-u_s^A(x)|\le2|A_i-A|
 =\frac{2|A_i^2-A^2|}{A_i+A}
 =\frac{4\sqrt2\,\norm r_2}{A_i+A}\le4\sqrt2\,\norm r_2,
\]
since $A>1$. Substituting this estimate into the maximum above and
taking the supremum over $x$, we obtain
\begin{equation}\label{eq:Bshape}
 \norm{v-u_s}_\infty\le4\sqrt2\,\norm r_2.
\end{equation}

We next estimate the first term $u-u_s$ in \eqref{eq:Bsplit}.
By \eqref{eq:Bprofile}, $u=u_0$ on $I$ and
$u_s(x)=u_0(x-s)$. Differentiating with respect to $s$ gives
$\partial_s u_s(x)=-u_0'(x-s)$. Integrating from $0$ to $s$, we obtain
\[
 u_s(x)-u_0(x)=-\int_0^s u_0'(x-\sigma)\dd\sigma.
\]
The $L^2(\R)$ norm of $u_0'(\cdot-\sigma)$ is independent of $\sigma$, so
\[
 \norm{u_s-u}_2
 \le\int_{\min(0,s)}^{\max(0,s)}\norm{u_0'(\cdot-\sigma)}_{L^2(\R)}\dd\sigma
 =|s|\norm{u_0'}_{L^2(\R)}.
\]
Using $u_0'(x)=-A^2(2\eps)^{-1}\operatorname{sech}^2(Ax/(2\eps))$
and the substitutions $z:=Ax/(2\eps)$ and $p:=\tanh z$, we find
\[
 \norm{u_0'}_{L^2(\R)}^2
 =\frac{A^3}{2\eps}\int_\R\operatorname{sech}^4 z\dd z
 =\frac{A^3}{2\eps}\int_{-1}^1(1-p^2)\dd p
 =\frac{2A^3}{3\eps},
\]
which gives
\begin{equation}\label{eq:BshiftQuadratic}
 \norm{u_s-u}_2^2\le\frac{2A^3}{3\eps}s^2.
\end{equation}

We now combine the bounds for the two terms in \eqref{eq:Bsplit}.
By the triangle inequality, we have
\[
 \norm{u-v}_2^2\le2\norm{u-u_s}_2^2+2\norm{u_s-v}_2^2.
\]
Since $|I|=2$, the estimate \eqref{eq:Bshape} gives
$2\norm{u_s-v}_2^2\le128\norm r_2^2$. Together with the estimate
\eqref{eq:BshiftQuadratic} for $\norm{u-u_s}_2^2$, this gives
\begin{equation}\label{eq:Btranslation}
 \norm{u-v}_2^2\le128\norm r_2^2+\frac{4A^3s^2}{3\eps}.
\end{equation}

To estimate $|s|$ in \eqref{eq:Btranslation}, we first consider $s\ge0$.
At the right endpoint, \eqref{eq:Bprofile} gives $u_s(1)=u(1-s)$, and
the boundary condition gives $v(1)=u(1)$. Hence
\[
 u_s(1)-v(1)=u(1-s)-u(1)=\int_{1-s}^{1}(-u'(x))\dd x.
\]
For $s<0$, we use the left endpoint, where $u_s(-1)=u(-1-s)$ and
$v(-1)=u(-1)$. In this case,
\[
 v(-1)-u_s(-1)=u(-1)-u(-1-s)=\int_{-1}^{-1-s}(-u'(x))\dd x.
\]
Both intervals of integration lie in $[-1,1]$ since $|s|<1$.
To estimate the integrals from below, we set $b:=-u'(1)$. By
\eqref{eq:Bsolution} and $A\tanh(A/(2\eps))=1$, we have, for $|x|\le1$,
\[
 -u'(x)=\frac{A^2}{2\eps}\operatorname{sech}^2\frac{Ax}{2\eps}
 \ge\frac{A^2}{2\eps}\operatorname{sech}^2\frac{A}{2\eps}
 =\frac{A^2-1}{2\eps}=b>0.
\]
Integrating the pointwise inequality $-u'(x)\ge b$ on the respective intervals gives
\[
 \begin{aligned}
 u_s(1)-v(1)&\ge\int_{1-s}^{1}b\dd x=bs &&(s\ge0),\\
 v(-1)-u_s(-1)&\ge\int_{-1}^{-1-s}b\dd x=b(-s)=b|s| &&(s<0).
 \end{aligned}
\]
Thus $b|s|\le\norm{u_s-v}_\infty$ in either case, and the estimate
\eqref{eq:Bshape} gives
\[
 |s|\le \frac{\norm{u_s-v}_\infty}{b}
 \le\frac{4\sqrt2}{b}\norm r_2.
\]

Substituting this bound for $|s|$ into \eqref{eq:Btranslation} and using
$A<2$, we obtain
\[
 \norm{v-u}_2^2\le128\norm r_2^2+\frac{4A^3s^2}{3\eps}
 \le C\left(1+\frac1{\eps b^2}\right)\norm r_2^2.
\]
To estimate $1/(\eps b^2)$, we first express $A-1$ in terms of $\eps$.
Using $\tanh z=(1-e^{-2z})/(1+e^{-2z})$ in $A\tanh(A/(2\eps))=1$ gives
\[
 A(1-e^{-A/\eps})=1+e^{-A/\eps},
 \qquad A-1=(A+1)e^{-A/\eps}.
\]
Since $1<A<2$ and $e^{1/\eps}\ge1/\eps$, it follows that
\[
 0<A-1\le3e^{-1/\eps},\qquad
 0<\frac{A-1}{\eps}\le\frac{3}{\eps}e^{-1/\eps}\le3.
\]
Substituting $A-1=(A+1)e^{-A/\eps}$ into $b=(A-1)(A+1)/(2\eps)$, we obtain
\[
 b=\frac{(A+1)^2}{2\eps}e^{-1/\eps}e^{-(A-1)/\eps}.
\]
The inequalities $2<A+1<3$ and $e^{-3}\le e^{-(A-1)/\eps}\le1$ now give
\[
 \frac{2e^{-3}}{\eps}e^{-1/\eps}\le b\le\frac{9}{2\eps}e^{-1/\eps},
 \qquad
 \frac{4}{81}\eps e^{2/\eps}\le\frac{1}{\eps b^2}
 \le\frac{e^6}{4}\eps e^{2/\eps}.
\]
Since $e^{2/\eps}\ge1+2/\eps$, we have $\eps e^{2/\eps}\ge2$ and hence
$1+1/(\eps b^2)\le C\eps e^{2/\eps}$ for every $0<\eps\le1$.
This proves the estimate \eqref{eq:Bstandard} when $\norm r_2\le1/(4\sqrt2)$.

We now consider the second case, $\norm r_2>1/(4\sqrt2)$.
To estimate $\norm{u-v}_2$, we first obtain an upper bound for $\norm v_\infty$.
By continuity, there is a point $x_*\in[-1,1]$ such that
$|v(x_*)|=\norm v_\infty$. If $\norm v_\infty>1$, the boundary values
$|v(\pm1)|=1$ imply $x_*\in I$ and hence $v'(x_*)=0$.

If $v(x_*)>1$, we use $v(-1)=1$ and choose
\[
 y:=\max\{x\in[-1,x_*]:v(x)=1\},\qquad v(x)>v(y)=1\quad(y<x\le x_*).
\]
If $v(x_*)<-1$, we use $v(1)=-1$ and choose
\[
 y:=\min\{x\in[x_*,1]:v(x)=-1\},\qquad v(x)<v(y)=-1\quad(x_*\le x<y).
\]
In each case the set defining $y$ is closed and contains the indicated
endpoint. For $x$ strictly between $y$ and $x_*$, we have
$(v(x)-v(y))/(x-y)>0$. Letting $x$ approach $y$ gives $v'(y)\ge0$.
By \eqref{eq:Bres}, $(-\eps v'+v^2/2)'=r$. Integrating from $y$ to $x_*$
and using $v'(x_*)=0$ and $(v(y))^2=1$, we obtain
\[
 \begin{aligned}
 \frac{\norm v_\infty^2}{2}=\frac{(v(x_*))^2}{2}
 &= -\eps v'(y)+\frac12+\int_y^{x_*}r(z)\dd z\\
 &\le \frac12+|x_*-y|^{1/2}\norm r_2
 \le \frac12+\sqrt2\,\norm r_2 .
 \end{aligned}
\]
Including the case $\norm v_\infty=1$, we conclude that
\begin{equation}\label{eq:Bamplitude}
 \norm v_\infty^2\le1+2\sqrt2\,\norm r_2.
\end{equation}
Since $|I|=2$, $\norm u_\infty=1$ and $\norm r_2>1/(4\sqrt2)$, it follows that
\[
 \norm{u-v}_2^2\le 2\norm u_2^2+2\norm v_2^2
 \le 8+8\sqrt2\,\norm r_2\le C\norm r_2^2.
\]
Using $\eps e^{2/\eps}\ge2$ once more gives \eqref{eq:Bstandard}
in the remaining case.

\smallskip
\noindent$\bullet$ \emph{On the sharpness of the $\eps e^{2/\eps}$ dependence.}
We fix $0<\eps\le1$. To prove that the factor in \eqref{eq:Bstandard} is
sharp, we seek an approximation for which the ratio of the squared error
to the standard loss is at least $c\eps e^{2/\eps}$.
For $v:=u+th$ with $h(-1)=h(1)=0$, the boundary values of $v$ are exact.
Substituting $v=u+th$ into \eqref{eq:Bres} and using \eqref{eq:Bproblem},
we obtain
\[
 r=t\bigl(-\eps h'+uh\bigr)'+t^2hh'.
\]
We choose $h$ so that the coefficient of $t$ is one. For this it suffices
to solve $-\eps h'+uh=x$. With $k:=A/(2\eps)$, the formula
\eqref{eq:Bsolution} gives $-u/\eps=2k\tanh(kx)$. Thus we can integrate
the equation for $h$ by writing
\[
 h'+2k\tanh(kx)h=-\frac{x}{\eps},
 \qquad
 \bigl(h(x)\cosh^2(kx)\bigr)'=-\frac{x}{\eps}\cosh^2(kx).
\]
Integrating from $x$ to $1$ and using $h(1)=0$, we obtain
\begin{equation}\label{eq:Bsharpdirection}
 h(x)=\frac{\operatorname{sech}^2(kx)}{\eps}\int_x^1 z\cosh^2(kz)\dd z
     =\frac{\operatorname{sech}^2(kx)}{\eps}\int_{|x|}^1 z\cosh^2(kz)\dd z,
\end{equation}
where the second equality follows from the oddness of $z\cosh^2(kz)$.

With this choice of $h$, we have $r=t+t^2hh'$. Since $|I|=2$,
\[
 \begin{aligned}
 \Lstd(v)=\norm r_2^2
 &=\int_I(t+t^2hh')^2\dd x\\
 &=2t^2+2t^3\underbrace{\int_I hh'\dd x}_{[h^2/2]_{-1}^{1}=0}
      +t^4\norm{hh'}_2^2\\
 &=t^2\bigl(2+t^2\norm{hh'}_2^2\bigr).
 \end{aligned}
\]
Using $v-u=th$ gives, for $t\ne0$,
\begin{equation}\label{eq:Bsharpratio}
 \frac{\norm{v-u}_2^2}{\Lstd(v)}
 =\frac{\norm h_2^2}{2+t^2\norm{hh'}_2^2}.
\end{equation}

We first estimate $\norm h_2^2$ from below. For $|x|\le\eps/2$, we have
$k|x|\le A/4<1/2$, and the interval $[1-\eps/2,1]$ is contained in
$[|x|,1]$. On this interval $z\ge1/2$ and
\[
 \cosh^2(kz)\ge\tfrac14 e^{2kz}
 \ge\tfrac14 e^{A/\eps-A/2}\ge\tfrac14 e^{1/\eps-1},
\]
where we used $1<A<2$. Restricting the integral in \eqref{eq:Bsharpdirection}
to this interval of length $\eps/2$, we obtain
\[
 h(x)\ge
 \frac{\operatorname{sech}^2(1/2)}{\eps}\,\frac{\eps}{2}\,
 \frac{e^{1/\eps-1}}{8}
 \ge c e^{1/\eps}\qquad (|x|\le\eps/2).
\]
Squaring and integrating over $[-\eps/2,\eps/2]$ gives
$\norm h_2^2\ge c\eps e^{2/\eps}$ for every $0<\eps\le1$.

For this fixed $\eps$, we choose $0<t_\eps\le e^{-1/\eps}$ small enough
that $t_\eps\norm{hh'}_2\le1$. For any $0<t\le t_\eps$, we set
$v_\eps:=u+th$. Since $t^2\norm{hh'}_2^2\le1$, the lower bound for
$\norm h_2^2$ and the identity \eqref{eq:Bsharpratio} give
\[
 c\eps e^{2/\eps}\le\frac{\norm h_2^2}{3}
 \le\frac{\norm{v_\eps-u}_2^2}{\Lstd(v_\eps)}
 \le C\eps e^{2/\eps},
\]
where the upper bound is \eqref{eq:Bstandard}. This proves sharpness.

\smallskip
\noindent$\bullet$ \emph{On the construction of $\tilde v$.}
We set $\tilde v:=u_{1/4}+g$, where $g$ is given by
\[
 g(x):=
 \begin{cases}
  \bigl(1-A\tanh\frac{5A}{8\eps}\bigr)(-2x-1)^3, & -1\le x<-\tfrac12,\\
  0, & -\tfrac12\le x\le\tfrac12,\\
  \bigl(A\tanh\frac{3A}{8\eps}-1\bigr)(2x-1)^3, & \tfrac12<x\le1.
 \end{cases}
\]
At $x=-1/2$ and $x=1/2$, the adjacent cubic and its first two derivatives
vanish, so $g\in C^2(\bar I)$. Furthermore,
\begin{align*}
 \tilde v(-1)
 &=u_{1/4}(-1)+g(-1)
  =A\tanh\frac{5A}{8\eps}+\left(1-A\tanh\frac{5A}{8\eps}\right)=1,\\
 \tilde v(1)
 &=u_{1/4}(1)+g(1)
  =-A\tanh\frac{3A}{8\eps}+\left(A\tanh\frac{3A}{8\eps}-1\right)=-1.
\end{align*}
Thus $\tilde v$ satisfies the boundary conditions exactly.

We next estimate the $L^2(I)$-norm of the PDE residual. The translate $u_{1/4}$
satisfies $-\eps u_{1/4}''+u_{1/4}u_{1/4}'=0$. Substituting
$\tilde v=u_{1/4}+g$ into the residual in \eqref{eq:Bres}, we therefore have
\[
 r(\tilde v)=-\eps g''+u_{1/4}g'+g u_{1/4}'+gg'.
\]
By the triangle inequality, $|u_{1/4}|\le A$, and $|I|=2$, we obtain
\begin{align*}
 \norm{r(\tilde v)}_2
 &\le\sqrt2\bigl(\eps\norm{g''}_\infty
       +(A+\norm g_\infty)\norm{g'}_\infty
       +\norm{g u_{1/4}'}_\infty\bigr)\\
 &\le C e^{-3/(4\eps)},
\end{align*}
where we used $A<2$, $0<\eps\le1$, and the bounds
\[
 \norm g_{W^{2,\infty}(I)}\le C e^{-3/(4\eps)},
 \qquad
 \sup_{x\in\operatorname{supp}g}|u_{1/4}'(x)|
 \le\frac{8}{\eps}e^{-1/(4\eps)}\le32.
\]
These follow by direct calculation from the formulas for $g$ and $u_{1/4}$,
using $|x-1/4|\ge1/4$ on $\operatorname{supp}g$ for the derivative bound.
Consequently, $\Lstd(\tilde v)\le C e^{-3/(2\eps)}$.

Finally, we obtain a lower bound for the error $\norm{\tilde v-u}_2$.
Restricting the integral to $[0,1/4]$ and applying the Cauchy--Schwarz
inequality, we have
\[
 \norm{\tilde v-u}_{L^2(I)}
 \ge\norm{\tilde v-u}_{L^2(0,1/4)}
 \ge2\left|\int_0^{1/4}(\tilde v-u)\dd x\right|.
\]
Since $g=0$ on $[0,1/4]$, we have
\[
 \tilde v(x)-u(x)
 =A\left(\tanh\frac{A(1/4-x)}{2\eps}+\tanh\frac{Ax}{2\eps}\right)
 \qquad(0\le x\le1/4).
\]
Using the substitution $x\mapsto1/4-x$ and $\cosh z\ge e^z/2$, we obtain
\begin{align*}
 \int_0^{1/4}(\tilde v-u)\dd x
 &=2A\int_0^{1/4}\tanh\frac{Ax}{2\eps}\dd x\\
 &=4\eps\log\cosh\frac{A}{8\eps}
 \ge\frac{A}{2}-4\eps\log2.
\end{align*}
Therefore,
\[
 \norm{\tilde v-u}_2
 \ge A-8\eps\log2\ge\tfrac12
 \qquad\left(0<\eps\le\frac{1}{16\log2}\right),
\]
where the last inequality uses $A>1$. This proves \eqref{eq:Bexample}.
\end{proof}

\begin{remark}\label{rem:Bsharpness}
As in the preceding two examples, the optimal constant in the standard
loss estimate~\eqref{eq:Bstandard} is large for small $\eps$. For each fixed
$0<\eps\le1$, the approximation $v_\eps=u+th$ in the proof attains the
optimal order $\eps e^{2/\eps}$ for the ratio of squared error to standard
loss. However, for $v_\eps$ we note that
\[
 \begin{aligned}
 \Lstd(v_\eps)&=t^2\bigl(2+t^2\norm{hh'}_2^2\bigr)
             \le3t^2=O(e^{-2/\eps}),\\
 \norm{v_\eps-u}_2^2&=t^2\norm h_2^2\le Ct=O(e^{-1/\eps}),
 \end{aligned}
\]
where we used $t\norm{hh'}_2\le1$, $t\le e^{-1/\eps}$, and the inequality
$\norm h_2^2\le C\norm{hh'}_2$, which follows by integration by parts since
$h(\pm1)=0$. Thus, this construction alone leaves
open whether the error can remain large at a very small standard loss.

By contrast, the approximation $\tilde v$ in~\eqref{eq:Bexample} has an
error of at least $1/2$ at an exponentially small standard loss for each
fixed sufficiently small $\eps>0$. We construct $\tilde v$ by shifting
the exact solution $u$ to the right by $1/4$ and adding $g$ near the
endpoints to restore the boundary values. The shifted solution still
satisfies the PDE, and its endpoint values differ from those of $u$ by
exponentially small amounts. The resulting correction and its first two
derivatives are therefore exponentially small, as is the PDE residual.
\end{remark}

In Theorem~\ref{prop:Bstandard}, shifting the exact solution and correcting
its values near the endpoints gives the approximation in~\eqref{eq:Bexample},
with a large error at a very small standard loss. It is therefore important
to control how far the layer in $v$ is shifted along the $x$-axis.
Indeed, if $u_s$ denotes the exact solution $u$ shifted by $s$, then
the size $|s|$ of the shift is controlled by the difference between
the integrals of $u_s$ and $u$ over $I$ through the inequality
\[
 |s|\le \left|\int_I u_s\dd x-\int_I u\dd x\right|,
 \qquad s\in I,\quad 0<\eps\le1.
\]
When $\norm r_2$ is sufficiently small and $s$ is the zero of $v$, the
estimate~\eqref{eq:Bshape} shows that $v$ is close to $u_s$ throughout
$I$, so their integrals are close as well. Training $v$ to match the
integral of $u$ is therefore a natural way to control the shift.

For problem~\eqref{eq:Bproblem}, we can determine $\int_I u\dd x$ from
the PDE and the boundary conditions alone. The equation
$\eps u''=u\,u'$ gives
$u'(x)=u'(-1)\exp\bigl(\eps^{-1}\int_{-1}^x u(y)\dd y\bigr)$.
Since $u(1)-u(-1)=-2$, we have $u'<0$ throughout $[-1,1]$.
Hence $u=\eps u''/u'=\eps(\log|u'|)'$. Moreover, integrating \eqref{eq:Bproblem}
shows that $-\eps u'+u^2/2$ is constant, and $(u(-1))^2=(u(1))^2=1$
then gives $u'(1)=u'(-1)$. Therefore
\begin{equation}\label{eq:Bmean}
 \int_Iu\dd x
 =\eps\log\frac{|u'(1)|}{|u'(-1)|}=0.
\end{equation}
Using \eqref{eq:Bmean}, we define the augmented loss $\Laug$ by
\begin{equation}\label{eq:Bloss}
 \Laug(v):=\norm r_2^2+|D(v)|^2,
\end{equation}
where
\[
 D(v):=\int_Iv\dd x.
\]

For this loss, we obtain the following estimate.
\begin{theorem}[Uniform estimate with an integral constraint]\label{prop:Baugmented}
For every $v\in H^2(I)$ with $v(-1)=1$ and $v(1)=-1$,
there is a constant $C$, independent of $v$ and $\eps$, such
that for every $0<\eps\le1$,
\begin{equation}\label{eq:Baugmented}
 \norm{v-u}_2^2\le C\min\left\{
       \eps^{-1}\Laug(v),\,\Laug(v)^{1/2}\right\}.
\end{equation}

\end{theorem}
\begin{proof}
We first consider $\norm r_2\le1/(4\sqrt2)$. For the unique zero $s$ of
$v$, we recall the estimate \eqref{eq:Btranslation}:
\begin{align}
 \norm{v-u}_2^2
 &\le2\norm{u-u_s}_2^2+2\norm{u_s-v}_2^2\label{eq:BaugmentedSplit}\\
 &\le128\norm r_2^2+\frac{4A^3s^2}{3\eps},\notag
\end{align}
where the last inequality uses the estimates \eqref{eq:Bshape} and
\eqref{eq:BshiftQuadratic} and $|I|=2$. In the proof of Theorem~\ref{prop:Bstandard},
$|s|$ was estimated
by $C\eps e^{1/\eps}\norm r_2$, with a factor that grows exponentially
as $\eps\to0$. Here we instead use $D(v)$, together with $\norm r_2$,
to improve the dependence of the error bound on $\eps$.

Since $u_s(x)=u_0(x-s)$ and $\int_I u_0\dd x=0$ by \eqref{eq:Bmean},
a change of variables gives, for every $s\in I$,
\begin{equation}\label{eq:Bintegralshift}
 \begin{aligned}
 \int_I u_s(x)\dd x
 &=\int_{-1-s}^{1-s}u_0(x)\dd x-\int_{-1}^{1}u_0(x)\dd x\\
 &=\int_{-1-s}^{-1}u_0(x)\dd x-\int_{1-s}^{1}u_0(x)\dd x.
 \end{aligned}
\end{equation}
For $0<s<1$, since $u_0$ is decreasing and $u_0(-1)=1$,
we have $u_0(x)\ge1$ on $[-1-s,-1]$. On $[1-s,1]$ we have $u_0(x)\le0$.
Thus \eqref{eq:Bintegralshift} gives
\[
 \int_I u_s(x)\dd x\ge\int_{-1-s}^{-1}1\dd x=s.
\]
For $s<0$, the identity $u_s(x)=-u_{-s}(-x)$ gives
$\int_I u_s\dd x=-\int_I u_{-s}\dd x\le s$. Hence
$|s|\le|\int_I u_s\dd x|$ for either sign of $s$, and equality holds
when $s=0$.

To express this bound in terms of $D(v)=\int_I v\dd x$, we use
\eqref{eq:Bshape} and $|I|=2$ to obtain
\[
 \left|\int_I u_s\dd x-D(v)\right|
 =\left|\int_I(u_s-v)\dd x\right|
 \le 2\norm{u_s-v}_\infty\le8\sqrt2\,\norm r_2.
\]
It follows that
\begin{equation}\label{eq:Bdisplacement}
 |s|\le\left|\int_I u_s\dd x\right|\le |D(v)|+8\sqrt2\,\norm r_2.
\end{equation}
Substituting \eqref{eq:Bdisplacement} for $|s|$ in \eqref{eq:Btranslation}
and using $A<2$ and $\eps\le1$, we have
\begin{align*}
 \norm{v-u}_2^2
 &\le 128\norm r_2^2+\frac{4A^3}{3\eps}
       \bigl(|D(v)|+8\sqrt2\,\norm r_2\bigr)^2\\
 &\le C\eps^{-1}\bigl(\norm r_2^2+|D(v)|^2\bigr)=C\eps^{-1}\Laug(v).
\end{align*}

Returning to the inequality \eqref{eq:BaugmentedSplit}, we can also obtain an
upper bound for $\norm{u-u_s}_2^2$ proportional to $|s|$, with a constant
independent of $\eps$.
Since $u_0$ is decreasing and
$u_s(x)=u_0(x-s)$, the difference $u_s-u$ has a constant sign on $I$.
Using the identities \eqref{eq:Bmean} and \eqref{eq:Bintegralshift},
together with $|u_0|\le A$, we have
\begin{align*}
 \int_I |u_s-u|\dd x
 &=\left|\int_I(u_s-u)\dd x\right|
  =\left|\int_I u_s\dd x\right|\\
 &\le\left|\int_{-1-s}^{-1}u_0(x)\dd x\right|
     +\left|\int_{1-s}^{1}u_0(x)\dd x\right|\\
 &\le2A|s|,
\end{align*}
which is illustrated in Figure~\ref{fig:burgers-decomposition}(b) for $s>0$.
Together with $|u_s-u|\le2A$, this gives
\begin{equation}\label{eq:BshiftLinear}
 \norm{u_s-u}_2^2\le 2A\int_I |u_s-u|\dd x\le4A^2|s|.
\end{equation}
We then use \eqref{eq:BshiftLinear} for $u-u_s$ and \eqref{eq:Bshape} for
$u_s-v$ in \eqref{eq:BaugmentedSplit}, together with the bound
\eqref{eq:Bdisplacement} for $|s|$, to obtain
\begin{align*}
 \norm{v-u}_2^2
 &\le 2\norm{u-u_s}_2^2+2\norm{u_s-v}_2^2\\
 &\le 8A^2|s|+128\norm r_2^2\\
 &\le 8A^2\bigl(|D(v)|+8\sqrt2\,\norm r_2\bigr)+128\norm r_2^2\\
 &\le C\bigl(|D(v)|+\norm r_2\bigr)\\
 &\le C\Laug(v)^{1/2},
\end{align*}
where we used $A<2$, $\norm r_2<1$, and the definition \eqref{eq:Bloss}.

It remains to consider $\norm r_2>1/(4\sqrt2)$. The estimate
\eqref{eq:Bamplitude}, together with $\norm u_\infty=1$ and $|I|=2$, gives
\begin{align*}
 \norm{v-u}_2^2
 &\le 2\norm v_2^2+2\norm u_2^2\\
 &\le 8+8\sqrt2\,\norm r_2\le C\norm r_2,
\end{align*}
where we used $1<4\sqrt2\,\norm r_2$. Since we also have
$\norm r_2\le4\sqrt2\,\norm r_2^2$, the definition \eqref{eq:Bloss}
and $0<\eps\le1$ yield
\begingroup
\postdisplaypenalty=10000
\[
 \norm{v-u}_2^2
 \le C\min\{\norm r_2^2,\norm r_2\}
 \le C\min\{\eps^{-1}\Laug(v),\Laug(v)^{1/2}\},
\]
which completes the proof.
\endgroup
\end{proof}

At a fixed loss value $\ell>0$, the inequality~\eqref{eq:Bstandard}
gives the upper bound $C\eps e^{2/\eps}\ell$ for the squared $L^2$ error,
while the inequality~\eqref{eq:Baugmented} gives the upper bound
$C\min\{\eps^{-1}\ell,\sqrt{\ell}\}$.
The latter is at most $C\sqrt{\ell}$ for every $0<\eps\le1$.
The growing difference between these bounds suggests a larger gain in accuracy at
smaller viscosity when the two losses reach comparable values.

\subsection{A reaction--diffusion equation with an unstable mean}\label{sec:reaction-diffusion}
The preceding examples concern time-independent equations. We now ask
whether an integral constraint can also improve error control in a
time-dependent problem and how its effect depends on the length of the
time interval.

For $T>0$, let $\Omega\subset\R^d$ be a bounded connected $C^2$ domain.
By rescaling space and the diffusion coefficient, we may assume $|\Omega|=1$.
Specifically, let us consider the following reaction--diffusion
equation for $u=u(t,x)$
\begin{equation}\label{eq:Dphysical}
 \begin{cases}
 u_t-\kappa\Delta u-\lambda u=f,&(t,x)\in(0,T)\times\Omega,\\
 \partial_nu=0,&(t,x)\in(0,T)\times\partial\Omega,\\
 u(0,x)=u_0(x),&x\in\Omega.
 \end{cases}
\end{equation}
Let $\mu_1$ denote the first positive Neumann eigenvalue of $-\Delta$. We assume $\kappa>0$ and $0<\lambda<\kappa\mu_1$.
The reaction amplifies the constant error mode at rate $\lambda$.
The condition $\lambda<\kappa\mu_1$ ensures that diffusion still damps
every mean-zero mode.
Further details on parabolic equations and related reaction--diffusion
models can be found in \citep{Brezis2011,Arcak2011}.

\leavevmode\llap{\enspace}For an approximation $v$ to $u$, we assume
$u,v\in H^1(0,T;L^2(\Omega))\cap L^2(0,T;H^2(\Omega))$,
with the Neumann condition satisfied exactly by both functions.
The initial condition is imposed through the loss.
Let $e:=v-u$. We define the residual $r$, the terms of the standard loss $\Lstd$, and the error norm $E$ by
\begin{align}\label{eq:Dloss}
 r&:=v_t-\kappa\Delta v-\lambda v-f,\notag\\
 \mathcal L_{\mathrm P}&:=\frac1T\int_0^T\|r\|_2^2\,dt,
 \qquad \mathcal L_{\mathrm{IC}}:=\|v(0)-u_0\|_2^2,\notag\\
 \Lstd&:=\mathcal L_{\mathrm P}+\mathcal L_{\mathrm{IC}},
 \qquad E^2:=\frac1T\int_0^T\|e\|_2^2\,dt.
\end{align}
The PDE residual loss and error are averaged over $(0,T)$.
The following theorem shows how the error bound for the standard loss depends on $T$.

\begin{theorem}[Standard loss estimate]\label{prop:Dstandard}
For every $T>0$ and every
$v\in H^1(0,T;L^2(\Omega))\cap L^2(0,T;H^2(\Omega))$
satisfying $\partial_nv=0$,
\begin{equation}\label{eq:Dstandard}
 E^2\le
 \max\left\{
 \frac{e^{2\lambda T}-1}{\lambda T},
 \frac{2(e^{\lambda T}-1)^2}{\lambda^2}
 \right\}\Lstd.
\end{equation}
Here, the exponential dependence $e^{2\lambda T}$ is sharp as
$T\to\infty$. For sufficiently large $T$, there is
$\widetilde v\in H^1(0,T;L^2(\Omega))\cap L^2(0,T;H^2(\Omega))$
satisfying $\partial_n\widetilde v=0$ and $\widetilde v(0)=u_0$ such that
\begin{equation}\label{eq:Dexample}
 \left(\frac1T\int_0^T\|\widetilde v-u\|_2^2\,dt\right)^{1/2}=1,
 \qquad \Lstd(\widetilde v)\le C e^{-2\lambda T},
\end{equation}
where $C$ depends only on $\lambda$ and is independent of $T$ and
$\widetilde v$.
\end{theorem}
\begin{proof}
$\bullet$ \emph{On the estimate~\eqref{eq:Dstandard}.}
By the definition of $r$ in~\eqref{eq:Dloss}, we have
$\partial_t v-\kappa\Delta v-\lambda v=f+r$. Subtracting the PDE for $u$
in~\eqref{eq:Dphysical} and using the Neumann boundary conditions, we obtain
\begin{equation}\label{eq:Derror}
 \partial_t e-\kappa\Delta e-\lambda e=r,\qquad \partial_n e=0.
\end{equation}
Multiplying the error equation~\eqref{eq:Derror} by $e$ and integrating
by parts over $\Omega$ gives, for almost every $t$,
\[
 \frac12\frac{d}{dt}\|e(t)\|_2^2
 =\lambda\|e(t)\|_2^2-\kappa\|\nabla e(t)\|_2^2
   +\int_\Omega r(t,x)e(t,x)\dd x.
\]
To estimate $\|e(t)\|_2$ even at times when it vanishes, we use
$q_\eta(t):=(\|e(t)\|_2^2+\eta^2)^{1/2}$ for $\eta>0$.
Since $q_\eta$ is absolutely continuous and $\|e(t)\|_2\le q_\eta(t)$,
the chain rule and the Cauchy--Schwarz inequality give
\[
 q_\eta'(t)
 \le \frac{\lambda\|e(t)\|_2^2+\|r(t)\|_2\|e(t)\|_2}{q_\eta(t)}
 \le \lambda q_\eta(t)+\|r(t)\|_2.
\]
It follows that $(e^{-\lambda t}q_\eta(t))'\le e^{-\lambda t}\|r(t)\|_2$.
Integrating from $0$ to $t$ and then letting $\eta\downarrow0$, we obtain
\[
 \|e(t)\|_2\le e^{\lambda t}\|e(0)\|_2
             +\int_0^t e^{\lambda(t-s)}\|r(s)\|_2\,ds.
\]
We square this inequality using $(a+b)^2\le 2a^2+2b^2$ and average over
$(0,T)$. Since $\|e(0)\|_2^2=\mathcal L_{\mathrm{IC}}$, integration of
$e^{2\lambda t}$ gives
\begin{equation}\label{eq:Daveraged}
 E^2\le\frac{e^{2\lambda T}-1}{\lambda T}\mathcal L_{\mathrm{IC}}
 +\frac2T\int_0^T
 \left(\int_0^t e^{\lambda(t-s)}\|r(s)\|_2\,ds\right)^2dt.
\end{equation}
For the integral in~\eqref{eq:Daveraged}, we first apply the
Cauchy--Schwarz inequality with weight $e^{\lambda(t-s)}$ to obtain
\[
 \left(\int_0^t e^{\lambda(t-s)}\|r(s)\|_2\,ds\right)^2
 \le \frac{e^{\lambda t}-1}{\lambda}
      \int_0^t e^{\lambda(t-s)}\|r(s)\|_2^2\,ds.
\]
Using $e^{\lambda t}\le e^{\lambda T}$ and changing the order of
integration on $0\le s\le t\le T$, we have
\[
 \begin{aligned}
 \frac1T\int_0^T\left(\int_0^t e^{\lambda(t-s)}\|r(s)\|_2\,ds\right)^2dt
 &\le \frac{e^{\lambda T}-1}{\lambda T}
      \int_0^T\|r(s)\|_2^2\int_s^T e^{\lambda(t-s)}\,dt\,ds\\
 &=\frac{e^{\lambda T}-1}{\lambda T}\int_0^T
      \frac{e^{\lambda(T-s)}-1}{\lambda}\|r(s)\|_2^2\,ds\\
 &\le \frac{(e^{\lambda T}-1)^2}{\lambda^2}\mathcal L_{\mathrm P}.
 \end{aligned}
\]
Substitution into~\eqref{eq:Daveraged} gives
\[
 E^2\le \frac{e^{2\lambda T}-1}{\lambda T}\mathcal L_{\mathrm{IC}}
       +\frac{2(e^{\lambda T}-1)^2}{\lambda^2}\mathcal L_{\mathrm P}.
\]
Since $\Lstd=\mathcal L_{\mathrm{IC}}+\mathcal L_{\mathrm P}$, taking
the larger of the two coefficients proves~\eqref{eq:Dstandard} for every $T>0$.

\smallskip
\noindent$\bullet$ \emph{On the existence of $\widetilde v$.}
We set $v(t,x):=u(t,x)+\sinh(\lambda t)$. Then $v(0)=u_0$,
$\partial_n v=0$, and the error equation~\eqref{eq:Derror} gives
\[
 r(t,x)=(\partial_t-\lambda)\sinh(\lambda t)
 =\lambda\bigl(\cosh(\lambda t)-\sinh(\lambda t)\bigr)
 =\lambda e^{-\lambda t}.
\]
Since $|\Omega|=1$, substitution into~\eqref{eq:Dloss} gives
\[
 \begin{aligned}
 E^2
 &=\frac1T\int_0^T\sinh^2(\lambda t)\,dt
   =\frac{\sinh(2\lambda T)-2\lambda T}{4\lambda T},\\
 \Lstd(v)
 &=\frac1T\int_0^T\lambda^2 e^{-2\lambda t}\,dt
   =\frac{\lambda(1-e^{-2\lambda T})}{2T}.
 \end{aligned}
\]
For fixed $\lambda>0$ and all sufficiently large $T$, we have
$2\lambda T+\tfrac12 e^{-2\lambda T}\le\tfrac14 e^{2\lambda T}$, so
$\sinh(2\lambda T)-2\lambda T\ge\tfrac14 e^{2\lambda T}$.
Dividing $E^2$ by $\Lstd(v)$, we therefore obtain
\[
 \frac{E^2}{\Lstd(v)}
 =\frac{\sinh(2\lambda T)-2\lambda T}{2\lambda^2(1-e^{-2\lambda T})}
 \ge\frac{e^{2\lambda T}}{8\lambda^2}.
\]
Since the coefficient in~\eqref{eq:Dstandard} is at most
$\max\{1,2/\lambda^2\}e^{2\lambda T}$ for $\lambda T\ge1$,
the exponential order in~\eqref{eq:Dstandard} is sharp.

To obtain the unit error in~\eqref{eq:Dexample}, we set
$\widetilde v:=u+(v-u)/E$. This function satisfies the initial and Neumann
boundary conditions, and for every sufficiently large $T$ we have
\[
 \frac1T\int_0^T\|\widetilde v-u\|_2^2\,dt
 =\frac{1}{TE^2}\int_0^T\|v-u\|_2^2\,dt=1,
 \qquad
 \Lstd(\widetilde v)=\frac{\Lstd(v)}{E^2}\le 8\lambda^2e^{-2\lambda T}.
 \qedhere
\]
\end{proof}

For $\delta>0$, the estimate~\eqref{eq:Dstandard} guarantees $E\le\delta$ whenever
\[
 \Lstd\le\delta^2
 \min\left\{
 \frac{\lambda T}{e^{2\lambda T}-1},
 \frac{\lambda^2}{2(e^{\lambda T}-1)^2}
 \right\}.
\]
In the construction for \eqref{eq:Dexample}, the error
is constant in space and grows in proportion to $\sinh(\lambda t)$ even though
the initial condition is exact. On a long interval, a prediction with
small standard loss may therefore reproduce the spatial variation of $u$
while its spatial mean drifts away from the exact mean.

To control this drift, we first determine the exact spatial mean from
the given data. We write $\bar z:=\int_\Omega z\dd x$ for the spatial
mean and $m(t):=\overline{u(t)}$ for that of the exact solution. Integrating \eqref{eq:Dphysical} over
$\Omega$ and using the Neumann boundary condition gives
$m'(t)-\lambda m=\bar f$. Let $m_0:=\overline{u_0}$ denote the spatial mean of the initial data. With $m(0)=m_0$, the solution is
\begin{equation}\label{eq:Dtarget}
 m(t)=e^{\lambda t}
 \left[m_0+\int_0^t e^{-\lambda s}\bar f(s)\,ds\right].
\end{equation}
We penalize deviations of $\bar v(t)$ from the mean in
\eqref{eq:Dtarget} by defining the augmented loss
\begin{equation}\label{eq:Dmoment}
 \Laug:=\Lstd+\frac1T\int_0^T|\bar v(t)-m(t)|^2\,dt.
\end{equation}
For this loss, we obtain the following estimate.

\begin{theorem}[Uniform estimate with an integral constraint]\label{prop:Daugmented}
For every $T>0$ and every
$v\in H^1(0,T;L^2(\Omega))\cap L^2(0,T;H^2(\Omega))$
satisfying $\partial_nv=0$,
\begin{equation}\label{eq:Duniform}
 E^2\le\max\left\{1,\frac{1}{(\kappa\mu_1-\lambda)^2}\right\}\Laug.
\end{equation}
\end{theorem}
\begin{proof}
We write $w:=e-\bar e$, where $\bar e(t)=\bar v(t)-m(t)$.
Since $|\Omega|=1$, we have $\int_\Omega w(t,x)\dd x=0$.
Expanding the square of $e(t)=w(t)+\bar e(t)$, we obtain
\[
 \begin{aligned}
 \|e(t)\|_2^2
 &=\|w(t)\|_2^2+2\bar e(t)\int_\Omega w(t,x)\dd x+|\bar e(t)|^2\\
 &=\|w(t)\|_2^2+|\bar e(t)|^2.
 \end{aligned}
\]
Averaging over $(0,T)$ gives
\begin{equation}\label{eq:DerrorSplit}
 E^2=\frac1T\int_0^T\|w(t)\|_2^2\,dt
     +\frac1T\int_0^T|\bar v(t)-m(t)|^2\,dt.
\end{equation}
The second term in~\eqref{eq:DerrorSplit} is already included
in the augmented loss~\eqref{eq:Dmoment}. To estimate the first term, we
integrate the error equation~\eqref{eq:Derror} over $\Omega$ to obtain
\[
 \bar e'(t)-\kappa\int_\Omega\Delta e(t,x)\dd x-\lambda\bar e(t)=\bar r(t).
\]
By the Neumann boundary condition, $\int_\Omega\Delta e\,\dd x
=\int_{\partial\Omega}\partial_n e\,\dd S=0$. Hence the spatial mean satisfies
\[
 \bar e'(t)-\lambda\bar e(t)=\bar r(t).
\]
Subtracting this equation from~\eqref{eq:Derror}, we obtain
\[
 \partial_t w-\kappa\Delta w-\lambda w=r-\bar r,\qquad \partial_n w=0.
\]
Multiplying by $w$ and integrating by parts, we have, for almost every $t$,
\[
 \frac12\frac{d}{dt}\|w(t)\|_2^2
 +\kappa\|\nabla w(t)\|_2^2-\lambda\|w(t)\|_2^2
 =\int_\Omega (r-\bar r)w\dd x.
\]
Since $\bar w=0$, the Poincar\'e inequality gives
$\|\nabla w(t)\|_2^2\ge\mu_1\|w(t)\|_2^2$. Using this on the left
and the Cauchy--Schwarz and Young inequalities on the right, we obtain
\[
 \begin{aligned}
 \frac12\frac{d}{dt}\|w(t)\|_2^2+(\kappa\mu_1-\lambda)\|w(t)\|_2^2
 &\le \|r(t)-\bar r(t)\|_2\|w(t)\|_2\\
 &\le \frac{\kappa\mu_1-\lambda}{2}\|w(t)\|_2^2
       +\frac{\|r(t)-\bar r(t)\|_2^2}{2(\kappa\mu_1-\lambda)},
 \end{aligned}
\]
where $\kappa\mu_1-\lambda>0$. Collecting the terms in $\|w(t)\|_2^2$ gives
\[
 \frac{d}{dt}\|w(t)\|_2^2+(\kappa\mu_1-\lambda)\|w(t)\|_2^2
 \le \frac{\|r(t)-\bar r(t)\|_2^2}{\kappa\mu_1-\lambda}.
\]
Multiplying by $e^{(\kappa\mu_1-\lambda)t}$ and integrating from $0$ to
$t$, we obtain
\[
 \|w(t)\|_2^2\le e^{-(\kappa\mu_1-\lambda)t}\|w(0)\|_2^2
 +\frac{1}{\kappa\mu_1-\lambda}
  \int_0^t e^{-(\kappa\mu_1-\lambda)(t-s)}\|r(s)-\bar r(s)\|_2^2\,ds.
\]
We now average over $(0,T)$. In the residual term the integrand is
nonnegative, so we may integrate first in $t$ from $s$ to $T$. This gives
\begin{equation}\label{eq:DwAverage}
 \begin{aligned}
 \frac1T\int_0^T\|w(t)\|_2^2\,dt
 &\le \frac{1-e^{-(\kappa\mu_1-\lambda)T}}{T(\kappa\mu_1-\lambda)}
       \|w(0)\|_2^2\\
 &\quad+\frac{1}{T(\kappa\mu_1-\lambda)^2}
       \int_0^T\bigl(1-e^{-(\kappa\mu_1-\lambda)(T-s)}\bigr)
                  \|r(s)-\bar r(s)\|_2^2\,ds.
 \end{aligned}
\end{equation}
The coefficient of $\|w(0)\|_2^2$ is at most one by $1-e^{-z}\le z$
for $z\ge0$, and $1-e^{-(\kappa\mu_1-\lambda)(T-s)}\le1$.
For $e(0)=w(0)+\bar e(0)$ and $r(s)=(r(s)-\bar r(s))+\bar r(s)$,
we expand the squares as in the derivation of~\eqref{eq:DerrorSplit}.
Since $\int_\Omega w(0,x)\dd x=0$ and $\int_\Omega(r(s,x)-\bar r(s))\dd x=0$,
we obtain
\[
 \begin{aligned}
 \|w(0)\|_2^2&=\|e(0)\|_2^2-|\bar e(0)|^2\le\mathcal L_{\mathrm{IC}},\\
 \frac1T\int_0^T\|r(s)-\bar r(s)\|_2^2\,ds
 &=\frac1T\int_0^T\bigl(\|r(s)\|_2^2-|\bar r(s)|^2\bigr)\,ds
 \le\mathcal L_{\mathrm P}.
 \end{aligned}
\]
Using the bounds for the initial error and the residual in~\eqref{eq:DwAverage}
and substituting the result into~\eqref{eq:DerrorSplit}, we obtain, for every $T>0$,
\[
 \begin{aligned}
 E^2&\le\mathcal L_{\mathrm{IC}}
       +\frac{\mathcal L_{\mathrm P}}{(\kappa\mu_1-\lambda)^2}
       +\frac1T\int_0^T|\bar v(t)-m(t)|^2\,dt\\
 &\le\max\left\{1,\frac{1}{(\kappa\mu_1-\lambda)^2}\right\}\Laug,
 \end{aligned}
\]
where the last inequality follows from the definition~\eqref{eq:Dmoment}
of $\Laug$.
\end{proof}

Penalizing the difference between $\bar v(t)$ and $m(t)$ in
\eqref{eq:Dmoment} gives the bound \eqref{eq:Duniform} on the
time-averaged error. Unlike the coefficient in \eqref{eq:Dstandard},
which grows exponentially with $T$, the bound \eqref{eq:Duniform}
is independent of $T$ when $\kappa\mu_1>\lambda$.
For the benchmark $\Omega=(0,1)$, $\kappa=\pi^{-2}$ and $\lambda=1/2$,
the first positive Neumann eigenvalue is $\mu_1=\pi^2$. Thus
$\kappa\mu_1-\lambda=1/2$, and the inequality~\eqref{eq:Duniform} gives the explicit
bound $E^2\le4\Laug$, independently of $T$.

The condition $\kappa\mu_1-\lambda>0$ is essential for controlling
the remaining mean-zero error.
Let $\phi_1$ be an $L^2(\Omega)$-normalized Neumann eigenfunction corresponding to $\mu_1$. At $\lambda=\kappa\mu_1$, choosing $e(t,x):=t\phi_1(x)$ gives $r=\phi_1$,
$\bar e=0$, $\mathcal L_{\mathrm{IC}}=0$, and
$E^2/\Lstd=T^2/3$. For $\lambda>\kappa\mu_1$, a mean-zero mode also grows in time.

\subsection{Common mechanism}\label{sec:mechanism}
In all four examples, the solution error can remain large even when the
standard loss is very small. The approximations constructed above show
what form this error can take in each problem. In the nonlinear Neumann
problem, the approximation can differ from the exact solution by a
constant shift. In the Helmholtz problem, an error in the $\phi_2$ mode
can persist. In both cases, an error of fixed size can give a standard
loss of order $\varepsilon^2$. In the Burgers problem, the transition
layer can be displaced, and in the reaction--diffusion problem, the
spatial mean of the approximation can drift away from the exact mean
over time. The constructions in \eqref{eq:Bexample} and
\eqref{eq:Dexample} exhibit these errors at exponentially small standard
loss for small viscosity and long time intervals, respectively.

Each integral constraint helps control the form of error that arises in
the corresponding problem. For a constant shift $c$ in the nonlinear
Neumann problem, $|D(u+c)|$ in \eqref{eq:neumann-aug} is comparable to
$|c|$. In the Helmholtz problem, $D(v)$ in \eqref{eq:helmholtz-aug}
measures the $\phi_2$ coefficient of the error. In the Burgers problem,
the inequality \eqref{eq:Bdisplacement} shows that $D(v)$ in
\eqref{eq:Bloss} and the residual together control the displacement of
the layer. In the reaction--diffusion problem, the added term in
\eqref{eq:Dmoment} is exactly the time average of the squared error in
the spatial mean. Each constraint follows from the PDE and the given
data, and adding it to the loss gives uniform control of errors of the
corresponding form. At comparable loss levels, this helps explain why the
gain in accuracy can be substantial when the standard loss controls such
errors only weakly and small when it already controls them well.

We also note that rescaling the standard loss does not remove this
insensitivity.
Multiplying the standard loss by a positive constant such as $\eps^{-2}$
scales all error components equally, so the relative sensitivity among
them is unchanged and the minimizers remain the same. The optimizers used
in Section~\ref{sec:numerics} can produce similar iterates and
approximation errors after such rescaling, despite the change in loss
magnitude.

\section{Numerical experiments}\label{sec:numerics}
Using stability estimates, Section~\ref{sec:theory} explained
when and why adding an integral constraint can improve prediction accuracy.
This section presents numerical experiments that support this interpretation.
Specifically, we compare the accuracy of
networks trained with each loss across a range of parameter values or
final times.
We also examine whether the error components identified as
poorly controlled by the standard loss indeed dominate the errors of
approximations trained with that loss.

The numerical settings of each experiment are given in
Appendix~\ref{app:numerical-setup}, and the code is available upon
request.

\subsection{A nonlinear elliptic equation with Neumann boundary conditions}\label{sec:Nnumerics}
Let us consider the equation \eqref{eq:N} on $\Omega=(0,1)$ with $g(s)=s+s^3$ and source
\[
 f_\eps(x)=\pi^2\cos(\pi x)
 +\eps\bigl[1+\cos(\pi x)+(1+\cos(\pi x))^3\bigr].
\]
Then the corresponding exact solution is $u_*(x)=1+\cos(\pi x)$.
With a neural network $N_\theta$, we construct the approximate solution
$v_\theta$ by $v_\theta:=1+2\tanh N_\theta$, which takes values in $(-1,3)$.
Thus the bounds~\eqref{eq:g} hold with $\alpha=1$ and $\beta=28$ on the common range
of the solution and approximations.
In this experiment, the boundary penalty is averaged over the two endpoints.
We define the loss functions as follows.
\begin{equation}\label{eq:pilotlosses}
 \begin{aligned}
 \mathcal L_{\mathrm{std},\eps}(v)&:=\int_0^1r^2\dd x+\tfrac12\bigl(|v'(0)|^2+|v'(1)|^2\bigr),\\
 \mathcal L_{\mathrm{aug},\eps}(v)&:=\mathcal L_{\mathrm{std},\eps}(v)+|D|^2,\\
 D&:=\int_0^1g(v)\dd x-\tfrac72,
 \end{aligned}
\end{equation}
where $r:=-v''+\eps g(v)-f_\eps$ is the residual of \eqref{eq:N}.
The integral target $7/2$ follows from \eqref{eq:neumann-identity}.

For each comparison in this section, we use five paired
initializations. The two networks in each pair start from identical
weights and biases.
Let $\mathcal L$ denote the loss used to train each network and let $E:=\norm{v_\theta-u_*}_{H^1(0,1)}$ denote its absolute $H^1$ error.
Table~\ref{tab:neumann-mean} reports the final relative $L^2$ error
$\norm{v_\theta-u_*}_{L^2(0,1)}/\norm{u_*}_{L^2(0,1)}$
and the final loss $\mathcal L$ as the mean $\pm$ sample standard deviation over
the five runs.
The last column gives the maximum of $E^2/\mathcal L$ over the recorded iterates
of all five runs, including the initial network. Errors and losses at the recorded
iterates are evaluated by quadrature rules independent of the training
points.
The same reporting convention is used in the following subsections.

\begin{table}[!htbp]
\centering\large\setlength{\tabcolsep}{4pt}
\caption{Results for the nonlinear elliptic equation with Neumann boundary conditions.}
\label{tab:neumann-mean}
\input{table/neumann_mean_table.tex}\end{table}

In Table~\ref{tab:neumann-mean}, the two methods show similar accuracy
at $\eps=1$. For smaller values of $\eps$, the mean relative $L^2$ error
grows progressively under standard training, while under augmented training
it stays close to its level at $\eps=1$. The stability estimates in
Section~\ref{sec:neumann} help explain this widening gap. At a given loss
level, the factor $\eps^{-2}$ in the standard estimate \eqref{eq:Nstandard}
allows a larger $H^1$ error for smaller $\eps$, whereas the bound for the
augmented loss in \eqref{eq:Nuniform} is independent of $\eps$.

\begin{figure}[!htbp]
\centering
\includegraphics[width=\linewidth]{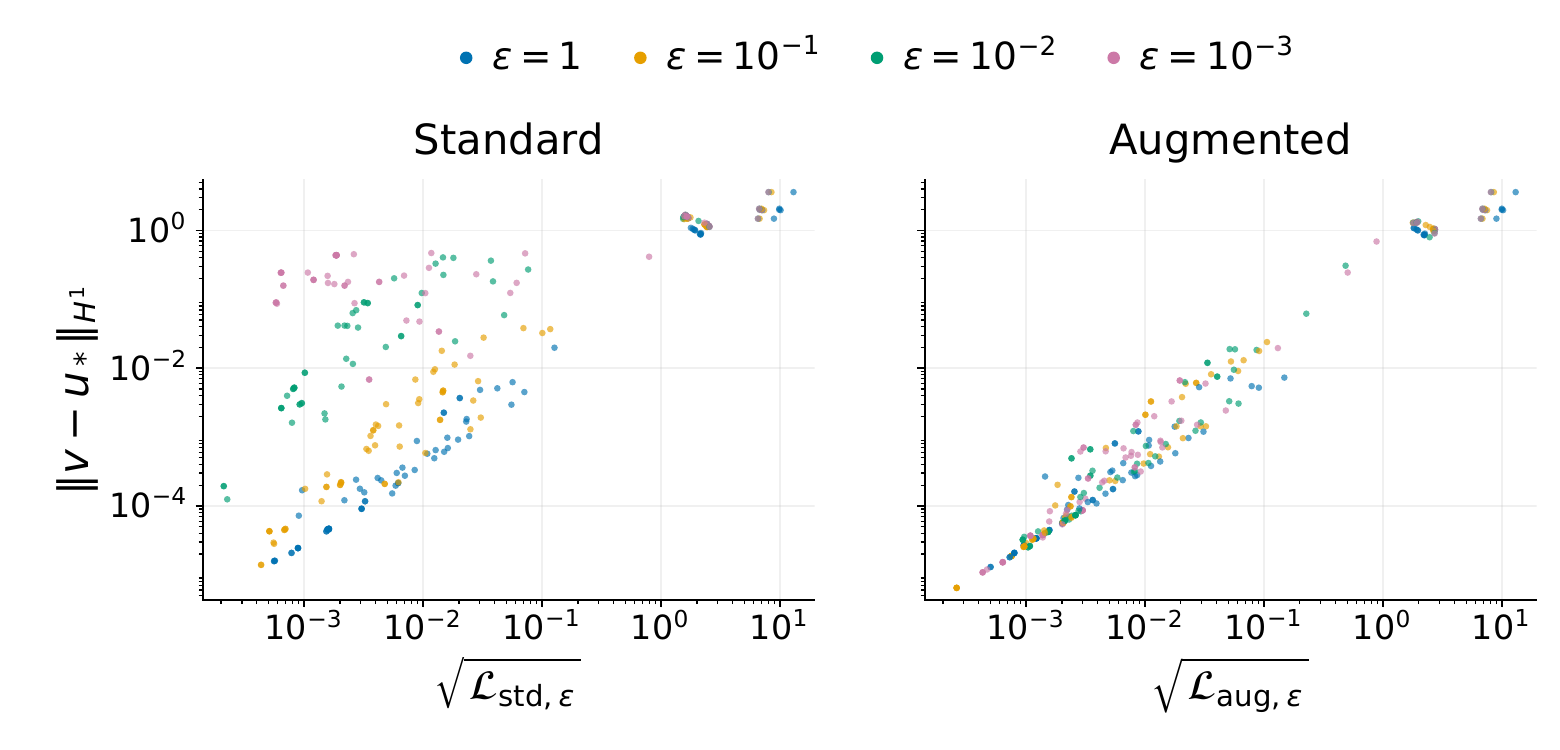}
\caption{Loss--error pairs for the nonlinear elliptic equation with Neumann boundary conditions.}
\label{fig:shared}
\end{figure}

Figure~\ref{fig:shared} plots the absolute $H^1$ error
appearing in the estimates \eqref{eq:Nstandard} and \eqref{eq:Nuniform}
against the square root of the loss used for training. Each point gives
the loss and error at one recorded iterate. At comparable small values
of the standard loss, the errors tend to be larger for smaller $\eps$.
With the augmented loss, the points for all four values of $\eps$ largely
overlap, and small losses are accompanied by small errors.

\begin{figure}[!htbp]
\centering
\includegraphics[width=\linewidth]{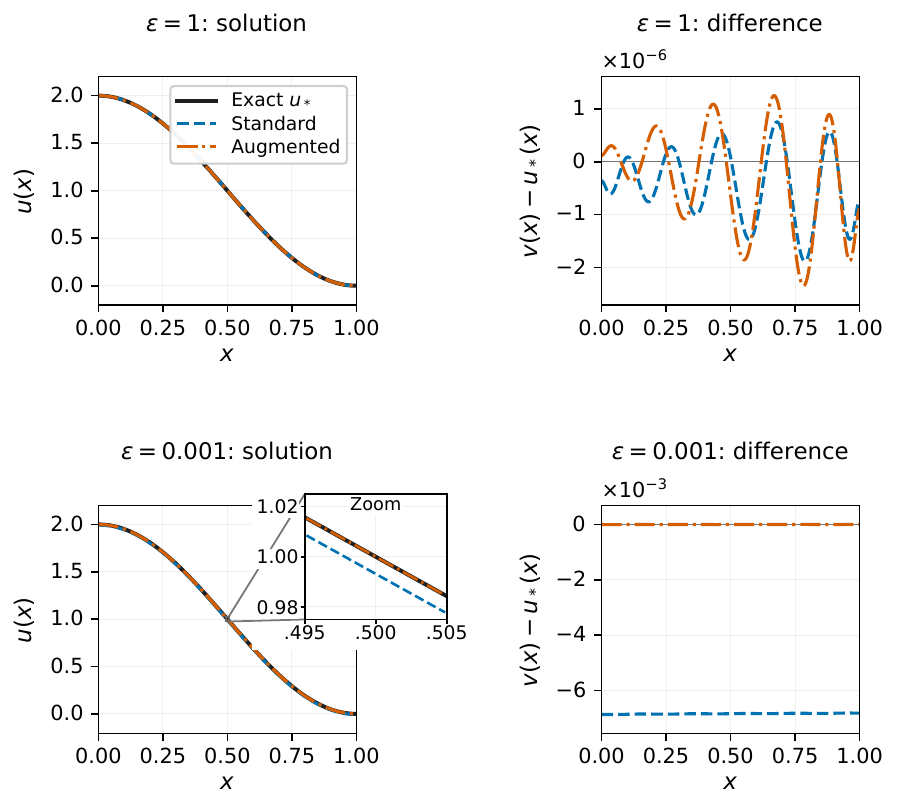}
\caption{Solution profiles from one pair of runs at each of $\eps=1$ and $\eps=10^{-3}$.}
\label{fig:neumann-profiles}
\end{figure}

Figure~\ref{fig:neumann-profiles} shows that both
approximations closely follow the exact solution at $\eps=1$.
At $\eps=10^{-3}$, the approximation trained with the standard loss lies
below the exact solution by nearly the same amount throughout the interval.
Adding the integral constraint substantially reduces this offset.
This contrast is consistent with the interpretation in
Section~\ref{sec:neumann} that the augmented loss provides stronger
control of constant shifts than the standard loss at small $\eps$.

\subsection{Helmholtz equation with Dirichlet boundary conditions}\label{sec:Helmnumerics}
Let us consider the equation \eqref{eq:HelmPDE} on $(0,1)$ with the fixed source
\[
 f(x)=-3\pi^2\sin(\pi x)+\frac{5\pi^2}{4}\sin(3\pi x).
\]
Then the corresponding exact solution is
\[
 u_*(x)=\frac{3\pi^2}{3\pi^2-\eps}\sin(\pi x)
 +\frac{5\pi^2}{4(5\pi^2+\eps)}\sin(3\pi x).
\]
With a neural network $N_\theta$, we construct the approximate solution
$v_\theta$ by $v_\theta(x):=x(1-x)N_\theta(2x-1)$, which imposes both
Dirichlet boundary conditions.
The standard and augmented losses are as defined
in \eqref{eq:HelmLoss} and \eqref{eq:helmholtz-aug}.
Since $\int_0^1 f(x)\phi_2(x)\,dx=0$, the integral identity in Section~\ref{sec:helmholtz}
gives $\int_0^1 u_*(x)\phi_2(x)\,dx=0$, so
$D(v)=\int_0^1v(x)\phi_2(x)\,dx$.
Here $E:=\norm{v_\theta-u_*}_{L^2(0,1)}$ denotes the absolute $L^2$ error.
The remaining reporting conventions follow Section~\ref{sec:Nnumerics}.

\begin{table}[!htbp]
\centering\large\setlength{\tabcolsep}{4pt}
\caption{Results for the Helmholtz equation with Dirichlet boundary conditions.}
\label{tab:helmholtz-mean}
\input{table/helmholtz_mean_table.tex}\end{table}

In Table~\ref{tab:helmholtz-mean}, the two methods show
similar accuracy at $\eps=1$, where the estimates \eqref{eq:HelmStandard}
and \eqref{eq:HelmUniform} give the same upper bound at a given loss value.
For smaller values of $\eps$, the two methods still attain comparable mean
final losses, but the mean relative $L^2$ error grows sharply under standard
training and remains small under augmented training. The largest recorded
$E^2/\mathcal L$ also grows sharply for the standard loss, while it stays
below one for the augmented loss. This contrast is consistent with the
bounds in Section~\ref{sec:helmholtz}, which give
$E^2/\mathcal L\leq\eps^{-2}$ for the standard loss and
$E^2/\mathcal L\leq 1$ for the augmented loss.

\begin{figure}[!htbp]
\centering
\includegraphics[width=\textwidth]{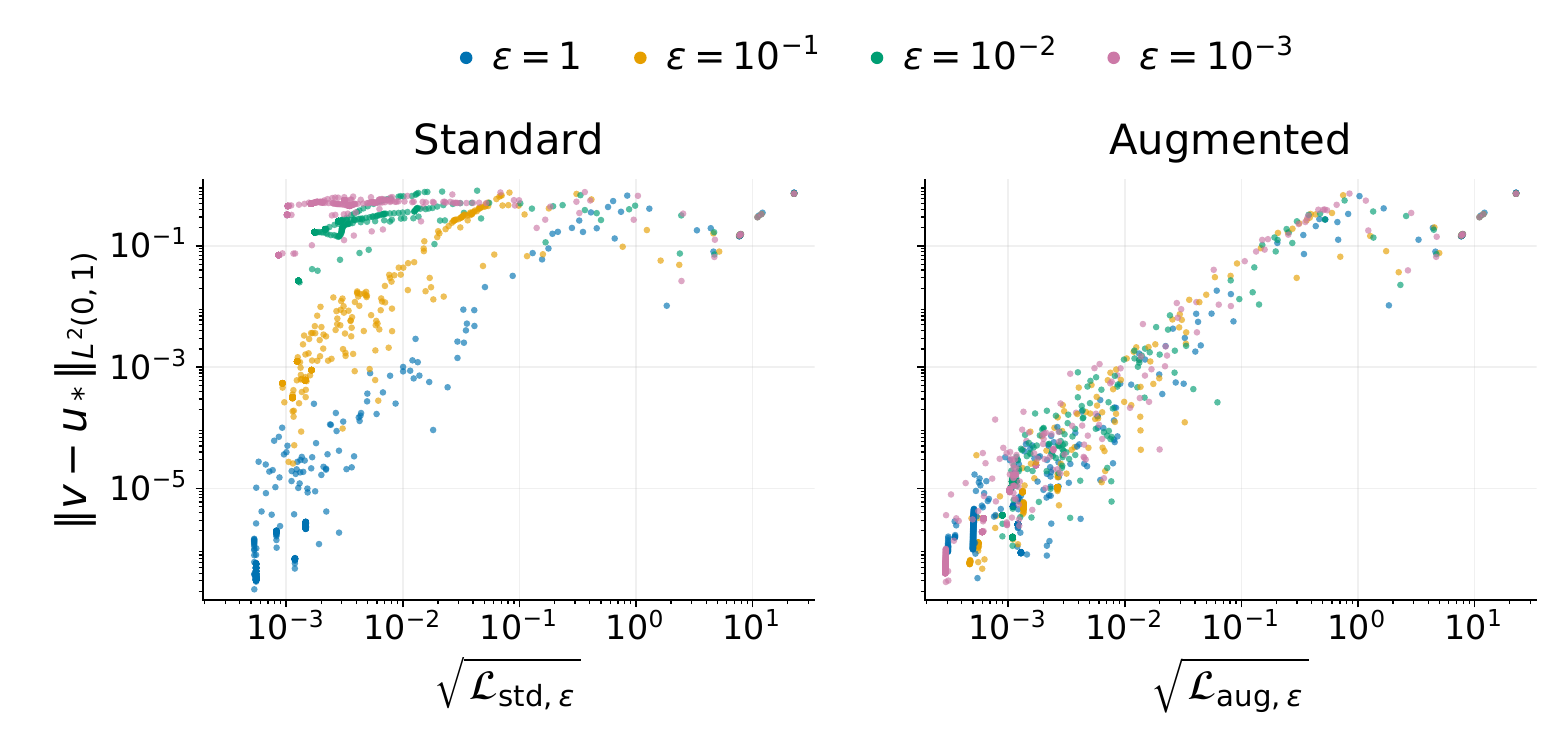}
\caption{Loss--error pairs for the Helmholtz equation with Dirichlet boundary conditions.}
\label{fig:helmholtz-candidates}
\end{figure}

Figure~\ref{fig:helmholtz-candidates} plots the absolute
$L^2$ error against the square root of the training loss. For the two
smallest values of $\eps$, many of the points obtained with the standard
loss form nearly horizontal bands, with substantial errors even at small
loss values. With the augmented loss, by contrast, smaller losses tend to
be accompanied by smaller errors for all four values of $\eps$, consistent
with the uniform estimate \eqref{eq:HelmUniform}.

\begin{figure}[!htbp]
\centering
\includegraphics[width=\linewidth]{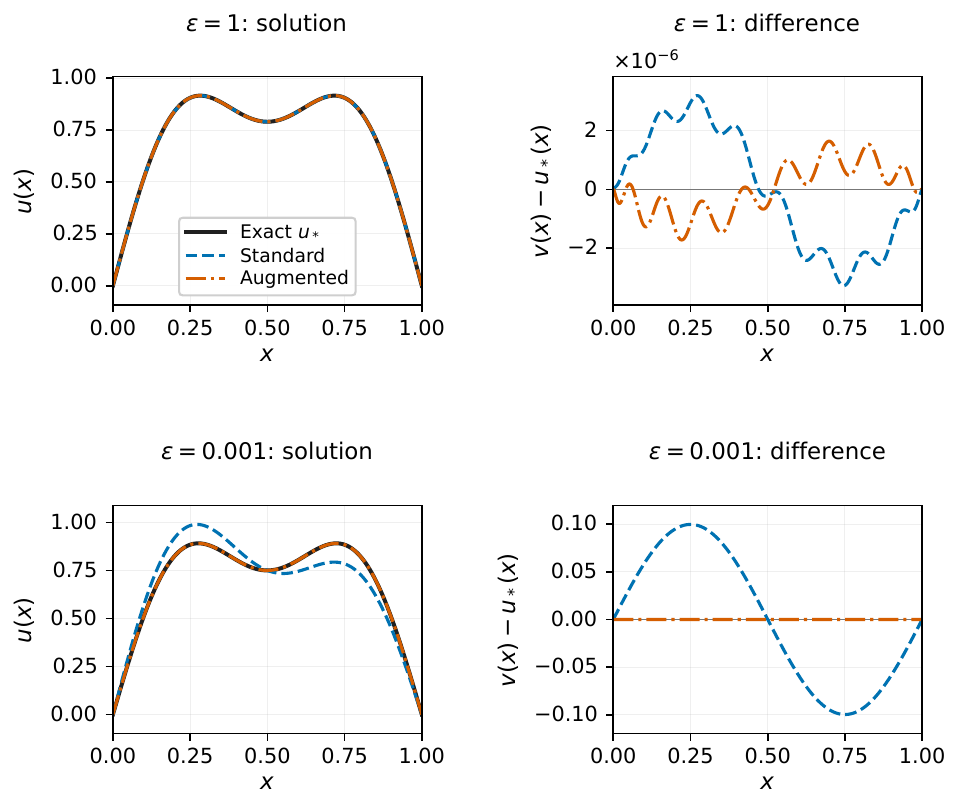}
\caption{Solution profiles from one pair of runs at each of $\eps=1$ and $\eps=10^{-3}$.}
\label{fig:helmholtz-profiles}
\end{figure}

Figure~\ref{fig:helmholtz-profiles} shows that both
approximations closely follow the exact solution at $\eps=1$.
At $\eps=10^{-3}$, the error under standard training is approximately a
multiple of $\sin(2\pi x)$. Adding the integral constraint substantially
reduces the amplitude of this oscillation. These profiles support the
interpretation in Section~\ref{sec:helmholtz} that the integral constraint
directly controls the sine component that the standard loss detects poorly
at small $\eps$.

\subsection[The stationary viscous Burgers equation]{The stationary viscous Burgers equation}\label{sec:Burgersnumerics}
In this subsection, we consider the equation \eqref{eq:Bproblem}.
Recall that the exact solution is
\[
 u_*(x)=-A\tanh\!\left(\frac{Ax}{2\eps}\right),
 \qquad A\tanh\!\left(\frac{A}{2\eps}\right)=1,\quad A>0.
\]
With a neural network $N_\theta$, we construct the approximate solution
$v_\theta$ by $v_\theta(x):=-x+(1-x^2)N_\theta(x)$, which imposes
both boundary conditions. The standard and augmented losses are as defined
in Section~\ref{sec:burgers},
\[
 \begin{aligned}
 \Lstd(v)&=\int_{-1}^1|-\eps v''+v\,v'|^2\,dx,\\
 \Laug(v)&=\Lstd(v)+|D(v)|^2,
 \qquad D(v)=\int_{-1}^1v\,dx.
 \end{aligned}
\]
The zero integral target follows from \eqref{eq:Bmean}.
Here $E:=\norm{v_\theta-u_*}_{L^2(I)}$ denotes the absolute $L^2$ error,
and the table reports both $E^2/\mathcal L$ and $E^4/\mathcal L$.
The remaining conventions follow Section~\ref{sec:Nnumerics}.

\begin{table}[!htbp]
\centering\large\setlength{\tabcolsep}{3pt}
\caption{Results for the stationary viscous Burgers equation.}
\label{tab:burgers-mean}
\newcommand{\pinnBurgersValue}[1]{#1}
\input{table/burgers_unit_weight_table.tex}\end{table}

\begin{samepage}
In Table~\ref{tab:burgers-mean}, the two methods give
similar mean relative $L^2$ errors at $\eps=1$. At the two smaller
viscosities, $\eps=0.1$ and $\eps=0.01$, the mean error grows sharply
under standard training and remains small under augmented training.
The largest recorded $E^2/\Lstd$ also increases as the viscosity decreases.
For the augmented loss, the largest recorded $E^4/\Laug$ remains below
one at all three viscosities, consistent with the viscosity-independent
bound on this ratio in \eqref{eq:Baugmented}.\par
\end{samepage}

\begin{figure}[!htbp]
\centering
\includegraphics[width=\textwidth]{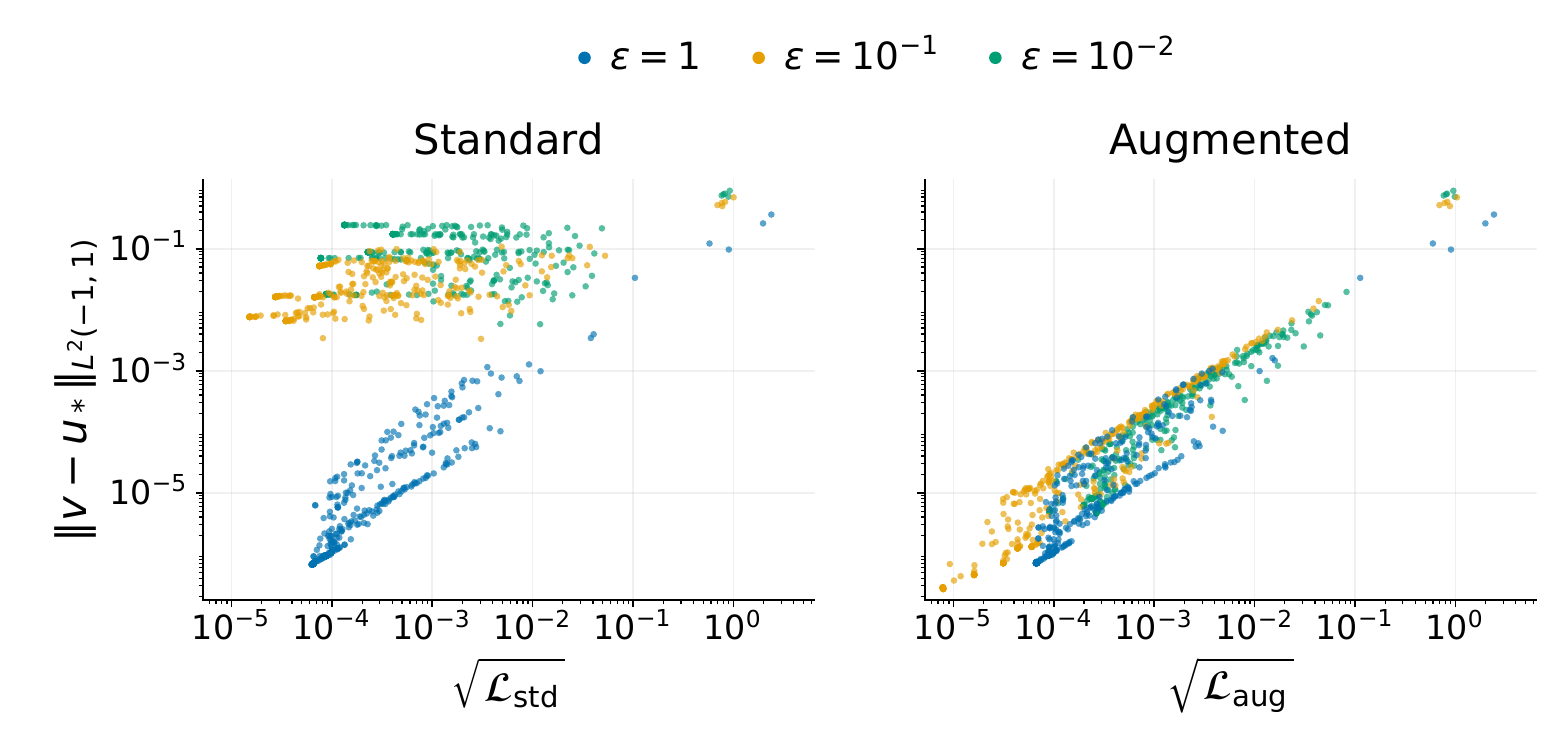}
\caption{Loss--error pairs for the stationary viscous Burgers equation.}
\label{fig:burgers-candidates}
\end{figure}

\begin{samepage}
Figure~\ref{fig:burgers-candidates} plots the absolute
$L^2$ error against the square root of the training loss. At comparable
small values of the standard loss, the errors at viscosities $\eps=0.1$
and $\eps=0.01$ are much larger than those at $\eps=1$. With the augmented
loss, by contrast, the error consistently decreases as the loss decreases,
regardless of the viscosity $\eps$, in agreement with the uniform estimate
\eqref{eq:Baugmented}.\par
\end{samepage}

\begin{figure}[!htbp]
\centering
\includegraphics[width=\linewidth]{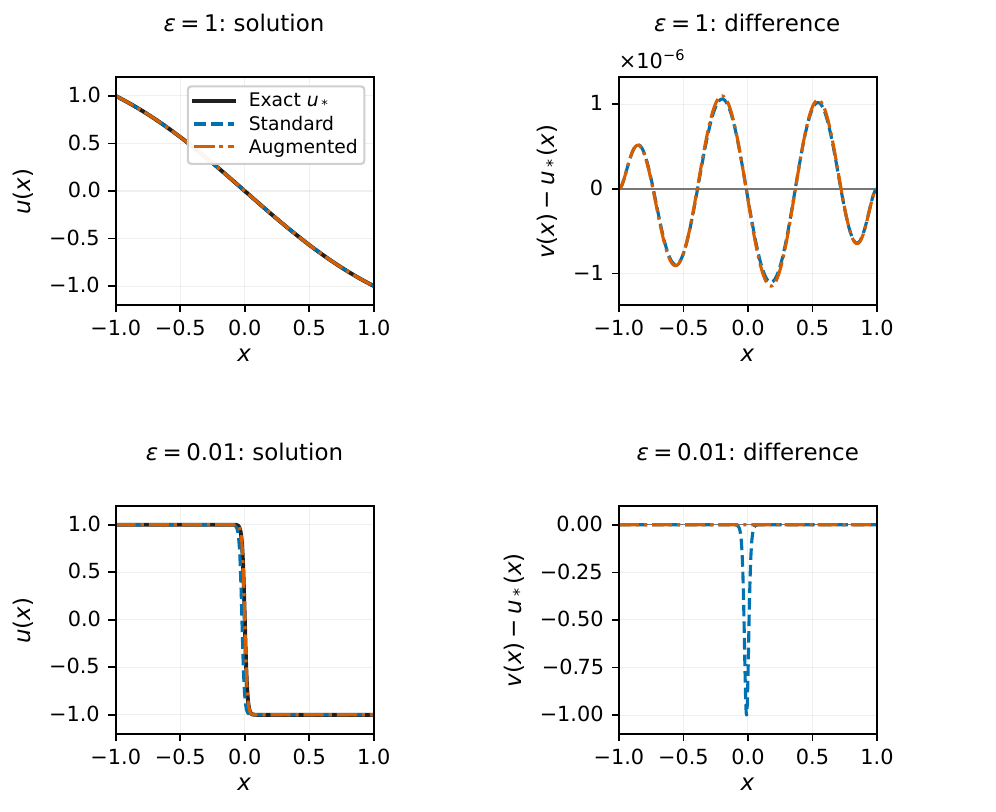}
\caption{Solution profiles from one pair of runs at each of $\eps=1$ and $\eps=10^{-2}$.}
\label{fig:burgers-profiles}
\end{figure}

\begin{samepage}
Figure~\ref{fig:burgers-profiles} shows that both
approximations closely follow the exact solution at viscosity $\eps=1$.
At $\eps=0.01$, the approximation trained with the standard loss captures
the sharp transition but places it slightly to the left of its exact
position at $x=0$. Adding the integral constraint brings the layer much
closer to $x=0$ and substantially reduces the error there. These profiles
support the interpretation in Section~\ref{sec:burgers} that, at small
viscosity $\eps$, a displaced layer can have a small standard loss and
the integral constraint improves control of its position.
\par
\end{samepage}

\subsection{A reaction--diffusion equation with Neumann boundary conditions}\label{sec:Dnumerics}
Let us consider the equation \eqref{eq:Dphysical} with $\Omega=(0,1)$,
$\kappa=\pi^{-2}$, $\lambda=1/2$ and
$u_0(x)=1+\tfrac12\cos(\pi x)$, with periodic source
\[
 f(t,x)=-\tfrac12+\tfrac12\cos t-\tfrac14\sin t
       +\bigl(\tfrac14\cos t-\tfrac12\sin t\bigr)\cos(\pi x).
\]
Then the corresponding exact solution is
\[
 u_*(t,x)=1+\tfrac12\sin t+\tfrac12\cos t\cos(\pi x).
\]
The coefficients, forcing and initial data are independent of $T$.
The spectral gap is $\kappa\mu_1-\lambda=1/2>0$.
The mean obtained from \eqref{eq:Dtarget} is
$m(t)=1+\tfrac12\sin t$.
This is a forced mean balance. Both the mean and the spatial variation
of the reference solution remain periodic as the terminal time increases.

With a neural network $N_\theta$, we construct the approximate solution
$v_\theta$ by $v_\theta(t,x):=N_\theta(2t/T-1,\cos(\pi x))$,
which satisfies the Neumann conditions exactly.
The standard and augmented losses are defined by
\begin{align}\label{eq:Dbenchmark-losses}
 \Lstd&:=\frac1T\int_0^T\int_0^1
       \left|v_t-\pi^{-2}v_{xx}-\tfrac12v-f(t,x)\right|^2\,dx\,dt
       \notag\\
    &\quad+\int_0^1\left|v(0,x)-1-\tfrac12\cos(\pi x)\right|^2\,dx,
       \notag\\
 \Laug&:=\Lstd+\frac1T\int_0^T\left|\int_0^1v(t,x)\,dx
                   -1-\tfrac12\sin t\right|^2\,dt.
\end{align}
With $e:=v_\theta-u_*$, we measure
\[
 E:=\left(\frac1T\int_0^T\|e(t)\|_2^2\,dt\right)^{1/2}
\]
and report the relative error
\[
 \frac{\|e\|_{L^2((0,T)\times(0,1))}}
      {\|u_*\|_{L^2((0,T)\times(0,1))}}.
\]
The added term in \eqref{eq:Dbenchmark-losses} is exactly
$T^{-1}\int_0^T|\bar e(t)|^2\,dt$.

\begin{table}[!htbp]
\centering\large\setlength{\tabcolsep}{4pt}
\caption{Results for the reaction--diffusion equation.}
\label{tab:Dmean}
\input{table/reaction_diffusion_error_table.tex}\end{table}

Table~\ref{tab:Dmean} shows that the two methods achieve similar accuracy
on the shortest interval, $T=4$. As the final time increases, the mean
relative error grows sharply under standard training and much more
modestly under augmented training. This widening gap is consistent with
the stability estimates in Section~\ref{sec:reaction-diffusion}. At a
given loss level, the error bound for the standard loss grows with $T$,
whereas the bound for the augmented loss is independent of $T$.

\begin{figure}[!htbp]
\centering
\includegraphics[width=\textwidth]{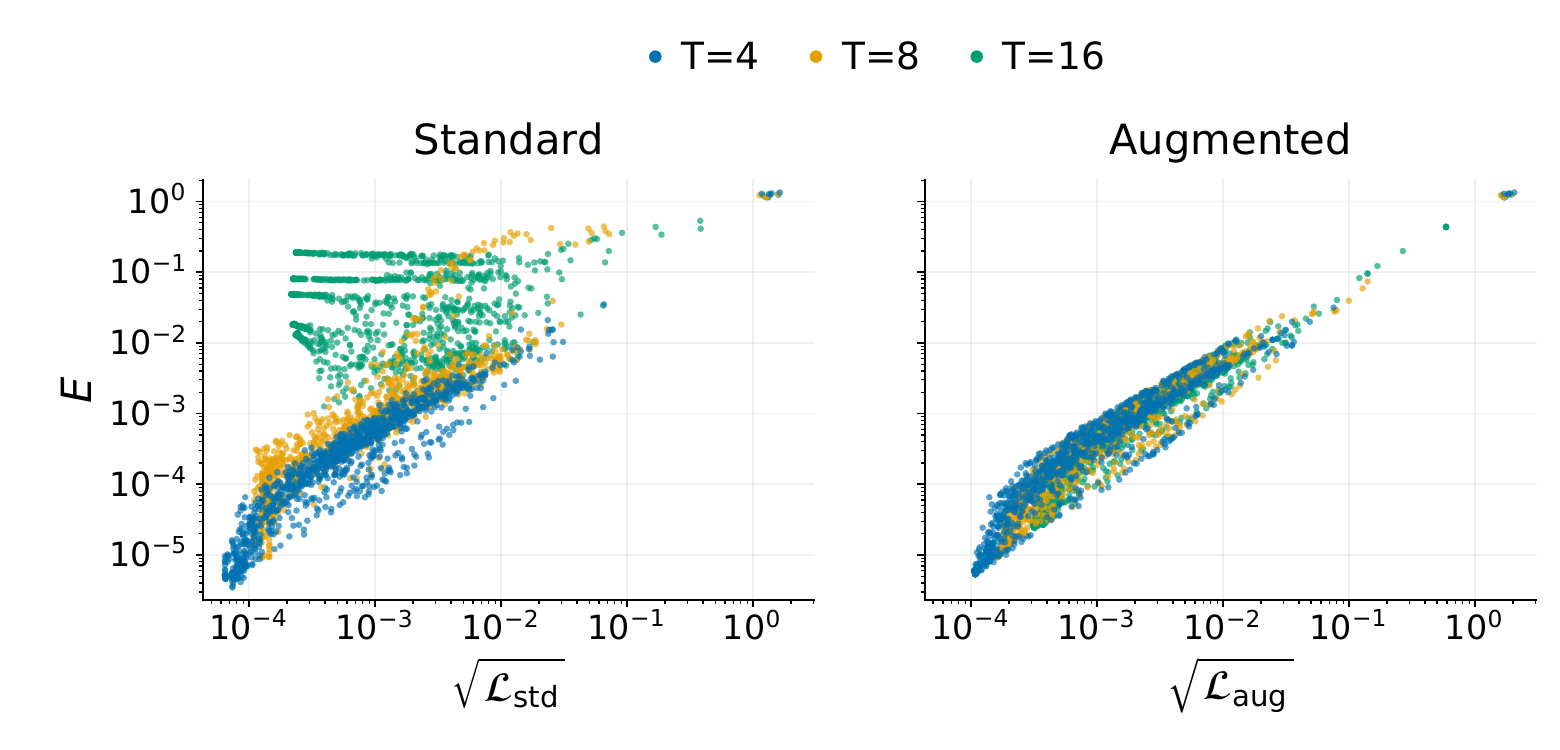}
\caption{Loss--error pairs for the reaction--diffusion equation.}
\label{fig:Dcandidates}
\end{figure}

Figure~\ref{fig:Dcandidates} plots the time-averaged $L^2$
error against the square root of the training loss. At comparably small
values of the standard loss, many of the errors at $T=16$ are substantially
larger than those at $T=4$, as the dependence on $T$ in the standard
estimate \eqref{eq:Dstandard} suggests. With the augmented loss, the
points for all three final times largely overlap, and smaller losses are
consistently accompanied by smaller errors, in agreement with the
$T$-independent estimate \eqref{eq:Duniform}.

\section{Conclusion}\label{sec:conclusion}
For each PDE, the stability estimates show that the
solution error can remain large even at small standard loss, depending on
the parameter or the final time, and identify a specific component that
dominates this error.
An integral constraint derived from the PDE and the given
data controls that component uniformly. This explains why adding the
constraint can improve accuracy in these regimes when training reduces
the two losses to similar values.

The experiments support this interpretation. The two methods achieve
similar accuracy when the stability estimates indicate that the standard
loss controls the solution error well. In regimes where the standard
loss provides poorer control of the solution error, however, training
with the augmented loss yields much smaller errors even when the two
methods attain similar loss values. Moreover, at small $\eps$, predictions by the standard loss
exhibit error shapes consistent with those predicted in
Section~\ref{sec:theory}. Adding the corresponding integral constraint
substantially reduces these errors.

We do not claim that every accuracy improvement from
adding integral constraints can be explained by this mechanism. Rather,
we demonstrate that stability theory provides such an explanation for
the representative cases studied here.
Extending this analysis to a broader range of PDEs is an important
direction for future work. It is equally important to analyze how
integral constraints affect the optimization landscape, training
convergence and generalization error, and to relate these effects to the
interpretation developed in this work.

\section*{Acknowledgement}
M.-S. Lee was supported by Basic Science Research Programs through the
National Research Foundation of Korea (NRF) funded by the Ministry of
Education (RS-2023-00244475 and RS-2024-00462755).

The author used GPT-6 Astra, Claude Opus 4.6, and
Gemini 4 Argon for English editing, and Claude Fable 5.2,
Claude Fable 5.1, and GPT-6 Astra to help write the code for the numerical
experiments. GPT-6 Astra also assisted with the construction of the
examples in Theorem~\ref{prop:Bstandard}. The author reviewed and rewrote
the corresponding proofs, independently verified all components of
the paper, and takes full responsibility for the content of this manuscript.

\appendix
\section{Numerical setup}\label{app:numerical-setup}
All computations use PyTorch in double precision on one NVIDIA RTX A6000
GPU per run. Table~\ref{tab:implementation} lists the settings for each
PDE. In the table, GL denotes Gauss--Legendre quadrature,
and \emph{Exact} indicates that the network satisfies the boundary
conditions by construction. All networks use $\tanh$
activations in the hidden layers.

\begin{table}[!t]
\centering\large\setlength{\tabcolsep}{3pt}
\caption{Numerical settings.}
\label{tab:implementation}
\input{table/training_settings.tex}\end{table}

For each PDE, both losses use the same optimizer settings. The Helmholtz
runs start directly with L-BFGS and use the full budget. In the other
three examples, Adam follows a cosine learning rate schedule. L-BFGS
then runs to the stated limit or until the infinity norm of the
parameter gradient is at most $10^{-10}$.

The Sobol counts in Table~\ref{tab:implementation} are per batch.
A new batch is drawn by randomized Sobol sampling every $100$ Adam
updates, and a fixed batch is used during L-BFGS.
For Helmholtz and Burgers, the PDE loss
and the integral constraint use the same fixed GL nodes and weights.
For Burgers, the $8192$ nodes come from an $8$-point GL rule on each
of $1024$ equal subintervals of $(-1,1)$.
For reaction--diffusion, the two quadrature orders for the integral
constraint refer to time and space, respectively.

\enlargethispage{2\baselineskip}
\end{document}

%% file: table/neumann_mean_table.tex
\renewcommand{\arraystretch}{1.15}
\begin{tabular*}{\linewidth}{@{\extracolsep{\fill}}ccrrr@{}}
\toprule
\multicolumn{1}{c}{$\varepsilon$} & \multicolumn{1}{c}{Loss} & \multicolumn{1}{c}{\shortstack{Relative\\$L^2$ error}} & \multicolumn{1}{c}{Final loss} & \multicolumn{1}{c}{$\max E^2/\mathcal L$} \\
\midrule
\multirow{2}{*}{1} & Standard & $(1.31\pm0.885)\times10^{-6}$ & $(1.31\pm1.03)\times10^{-6}$ & $3.76\times10^{-1}$ \\
 & Augmented & $(1.22\pm1.03)\times10^{-6}$ & $(1.45\pm1.53)\times10^{-6}$ & $3.60\times10^{-1}$ \\
\midrule
\multirow{2}{*}{0.1} & Standard & $(7.42\pm7.14)\times10^{-5}$ & $(1.46\pm1.70)\times10^{-6}$ & $1.55\times10^{0}$ \\
 & Augmented & $(1.06\pm0.852)\times10^{-6}$ & $(1.40\pm1.62)\times10^{-6}$ & $5.28\times10^{-1}$ \\
\midrule
\multirow{2}{*}{0.01} & Standard & $(3.16\pm2.56)\times10^{-3}$ & $(5.97\pm3.84)\times10^{-7}$ & $1.25\times10^{3}$ \\
 & Augmented & $(2.14\pm0.936)\times10^{-6}$ & $(3.10\pm2.48)\times10^{-6}$ & $5.20\times10^{-1}$ \\
\midrule
\multirow{2}{*}{0.001} & Standard & $(1.52\pm1.34)\times10^{-1}$ & $(3.39\pm5.20)\times10^{-6}$ & $1.45\times10^{5}$ \\
 & Augmented & $(2.03\pm2.58)\times10^{-6}$ & $(2.49\pm3.60)\times10^{-6}$ & $6.12\times10^{-1}$ \\
\bottomrule
\end{tabular*}

%% file: table/helmholtz_mean_table.tex
\renewcommand{\arraystretch}{1.15}
\begin{tabular*}{\linewidth}{@{\extracolsep{\fill}}ccrrr@{}}
\toprule
\multicolumn{1}{c}{$\varepsilon$} & \multicolumn{1}{c}{Loss} & \multicolumn{1}{c}{\shortstack{Relative\\$L^2$ error}} & \multicolumn{1}{c}{Final loss} & \multicolumn{1}{c}{$\max E^2/\mathcal L$} \\
\midrule
\multirow{2}{*}{1} & Standard & $(1.93\pm1.27)\times10^{-6}$ & $(9.67\pm8.07)\times10^{-7}$ & $8.88\times10^{-1}$ \\
 & Augmented & $(2.76\pm2.41)\times10^{-6}$ & $(9.36\pm7.19)\times10^{-7}$ & $4.16\times10^{-1}$ \\
\midrule
\multirow{2}{*}{0.1} & Standard & $(9.88\pm4.91)\times10^{-4}$ & $(1.71\pm0.728)\times10^{-6}$ & $9.71\times10^{1}$ \\
 & Augmented & $(7.17\pm5.97)\times10^{-6}$ & $(2.22\pm2.78)\times10^{-6}$ & $8.48\times10^{-1}$ \\
\midrule
\multirow{2}{*}{0.01} & Standard & $(2.45\pm1.28)\times10^{-1}$ & $(5.03\pm2.83)\times10^{-6}$ & $9.28\times10^{3}$ \\
 & Augmented & $(1.13\pm0.731)\times10^{-5}$ & $(1.13\pm0.209)\times10^{-6}$ & $8.08\times10^{-1}$ \\
\midrule
\multirow{2}{*}{0.001} & Standard & $(5.26\pm2.68)\times10^{-1}$ & $(3.09\pm3.91)\times10^{-6}$ & $1.94\times10^{5}$ \\
 & Augmented & $(1.07\pm1.33)\times10^{-5}$ & $(9.72\pm11.8)\times10^{-7}$ & $7.51\times10^{-1}$ \\
\bottomrule
\end{tabular*}

%% file: table/burgers_unit_weight_table.tex
\begingroup
\renewcommand{\arraystretch}{1.15}
\begin{tabular*}{\linewidth}{@{\extracolsep{\fill}}ccrrrr@{}}
\toprule
\multicolumn{1}{c}{$\varepsilon$} & \multicolumn{1}{c}{Loss} & \multicolumn{1}{c}{\shortstack{Relative\\$L^2$ error}} & \multicolumn{1}{c}{Final loss} & \multicolumn{1}{c}{$\max E^2/\mathcal L$} & \multicolumn{1}{c}{$\max E^4/\mathcal L$} \\
\midrule
\multirow{2}{*}{1} & Standard & $(9.12\pm1.91)\times10^{-7}$ & $(5.95\pm2.53)\times10^{-9}$ & $1.17\times10^{-1}$ & -- \\
 & Augmented & $(9.37\pm1.83)\times10^{-7}$ & $(6.24\pm2.42)\times10^{-9}$ & $9.53\times10^{-2}$ & $2.94\times10^{-3}$ \\
\midrule
\multirow{2}{*}{0.1} & Standard & $(1.56\pm1.47)\times10^{-2}$ & $(2.46\pm2.44)\times10^{-9}$ & $4.79\times10^{5}$ & -- \\
 & Augmented & \pinnBurgersValue{$(6.30\pm3.67)\times10^{-7}$} & \pinnBurgersValue{$(1.36\pm1.45)\times10^{-9}$} & \pinnBurgersValue{$5.55\times10^{-1}$} & \pinnBurgersValue{$2.15\times10^{-1}$} \\
\midrule
\multirow{2}{*}{0.01} & Standard & $(8.46\pm6.41)\times10^{-2}$ & $(4.94\pm6.64)\times10^{-8}$ & $3.40\times10^{6}$ & -- \\
 & Augmented & \pinnBurgersValue{$(4.30\pm0.963)\times10^{-6}$} & \pinnBurgersValue{$(5.13\pm2.91)\times10^{-8}$} & \pinnBurgersValue{$9.41\times10^{-1}$} & \pinnBurgersValue{$6.87\times10^{-1}$} \\
\bottomrule
\end{tabular*}
\endgroup

%% file: table/reaction_diffusion_error_table.tex
\renewcommand{\arraystretch}{1.15}
\begin{tabular*}{\linewidth}{@{\extracolsep{\fill}}ccrrr@{}}
\toprule
\multicolumn{1}{c}{$T$} & \multicolumn{1}{c}{Loss} & \multicolumn{1}{c}{\shortstack{Relative\\$L^2$ error}} & \multicolumn{1}{c}{Final loss} & \multicolumn{1}{c}{$\max E^2/\mathcal L$} \\
\midrule
\multirow{2}{*}{4} & Standard & $(4.99\pm1.02)\times10^{-6}$ & $(6.33\pm1.62)\times10^{-9}$ & $1.50\times10^{0}$ \\
 & Augmented & $(5.51\pm0.964)\times10^{-6}$ & $(1.48\pm0.344)\times10^{-8}$ & $7.97\times10^{-1}$ \\
\midrule
\multirow{2}{*}{8} & Standard & $(6.60\pm3.74)\times10^{-5}$ & $(1.65\pm0.320)\times10^{-8}$ & $1.25\times10^{3}$ \\
 & Augmented & $(1.05\pm0.160)\times10^{-5}$ & $(3.09\pm0.393)\times10^{-8}$ & $7.69\times10^{-1}$ \\
\midrule
\multirow{2}{*}{16} & Standard & $(6.15\pm6.29)\times10^{-2}$ & $(5.36\pm0.443)\times10^{-8}$ & $6.49\times10^{5}$ \\
 & Augmented & $(2.53\pm0.354)\times10^{-5}$ & $(1.25\pm0.240)\times10^{-7}$ & $7.72\times10^{-1}$ \\
\bottomrule
\end{tabular*}

%% file: table/training_settings.tex
\renewcommand{\arraystretch}{1.25}
\begin{tabular*}{\linewidth}{@{\extracolsep{\fill}}lcccc@{}}
\toprule
\textbf{Setting} & \textbf{Neumann} & \textbf{Helmholtz} & \textbf{Burgers} & \textbf{Reaction--diffusion} \\
\midrule
Hidden layers $\times$ width & $4\times64$ & $3\times48$ & $4\times64$ & $4\times64$ \\
\raisebox{0.5em}{Initialization} &
\shortstack{Xavier\\(Uniform)} &
\shortstack{Xavier\\(Normal)} &
\shortstack{Xavier\\(Uniform)} &
\shortstack{Xavier\\(Uniform)} \\
\midrule[0.35pt]
Adam updates & $10\,000$ & --- & $50\,000$ & $200\,000$ \\
Adam learning rate & $10^{-3}\to10^{-4}$ & --- & $10^{-3}\to10^{-5}$ & $10^{-3}\to10^{-8}$ \\
L-BFGS budget & $\le 10\,000$ & $2\,000$ & $\le 10\,000$ & $\le 10\,000$ \\
\midrule[0.35pt]
PDE loss points & $512$ Sobol & $96$ GL & $8192$ GL & $1024$ Sobol \\
Boundary loss points & $2$ endpoints & Exact & Exact & Exact \\
Initial condition points & --- & --- & --- & $64$ GL \\
Integral constraint points & $128$ GL & $96$ GL & $8192$ GL & $128\times32$ GL \\
\bottomrule
\end{tabular*}